\documentclass[11pt,a4paper]{article}

\usepackage[a4paper,left=1.65cm,right=1.65cm,top=1.65cm,bottom=1.65cm]{geometry}
\usepackage{t1enc}
\usepackage[utf8]{inputenc}
\usepackage{amsthm,amssymb}

\usepackage[mathcal,mathscr]{eucal}
\usepackage{enumerate}
\usepackage[sort,nocompress]{cite}
\usepackage{indentfirst}
\usepackage{todonotes}
\usepackage{subcaption}
\usepackage{hyperref}
\hypersetup{
    colorlinks=true,       % false: boxed links; true: colored links
    linkcolor=blue,        % color of internal links (change box color with linkbordercolor)
    citecolor=red,         % color of links to bibliography
    filecolor=magenta,     % color of file links
    urlcolor=cyan,         % color of external links
    linktocpage=true
}

\theoremstyle{plain}
\newtheorem{theorem}{Theorem}
\newtheorem{lemma}[theorem]{Lemma}
\newtheorem{corollary}[theorem]{Corollary}
\newtheorem{claim}[theorem]{Claim}

\theoremstyle{definition}

\newtheorem{remark}[theorem]{Remark}

\def\ol{\overline}

\usepackage{amsmath}

\begin{document}

\title{Supermodularity in Unweighted Graph Optimization III:  \break
Highly-connected Digraphs }

\author{Krist\'of B\'erczi\thanks{ MTA-ELTE Egerv\'ary Research Group,
Department of Operations Research, E\"otv\"os University, P\'azm\'any
P. s. 1/c, Budapest, Hungary, H-1117. e-mail:  {\tt berkri\char'100
cs.elte.hu .}} \ \ and  \ {Andr\'as Frank\thanks{MTA-ELTE Egerv\'ary
Research
Group, Department of Operations Research, E\"otv\"os University,
P\'azm\'any P. s. 1/c, Budapest, Hungary, H-1117. e-mail:  {\tt
frank\char'100 cs.elte.hu .} } } }

\maketitle

\begin{abstract} By generalizing a recent result of Hong, Liu, and Lai
\cite{Hong-Liu-Lai} on characterizing the degree-sequences of simple
strongly connected directed graphs, a characterization is provided for
degree-sequences of simple $k$-node-connected digraphs.  More
generally, we solve the directed node-connectivity augmentation
problem when the augmented digraph is degree-specified and simple.  As
for edge-connectivity augmentation, we solve the special case when the
edge-connectivity is to be increased by one and the augmenting digraph
must be simple.  \end{abstract}

\section{Introduction} \label{intro}

There is an extensive literature of problems concerning degree
sequences of graphs or digraphs with some prescribed properties such
as simplicity or $k$-connectivity.  For example, Edmonds
\cite{Edmonds64} characterized the degree-sequences of simple
$k$-edge-connected undirected graphs, while Wang and Kleitman
\cite{Wang-Kleitman} solved the corresponding problem for simple
$k$-node-connected graphs.

In what follows, we consider only directed graphs for which the
default understanding will be throughout the paper that loops and
parallel arcs are allowed.  When neither loops nor parallel arcs from
$u$ to $v$ are allowed we speak of {\bf simple} digraphs.  Oppositely
oriented arcs $uv$ and $vu$, however, are allowed in simple digraphs.
A typical problem is as follows.  Given an $n$-element ground-set $V$,
decide for a specified integer-valued function $m_i:V\rightarrow \mathbb{Z}_+$ if there is a digraph $D=(V,A)$ with some prescribed properties
realizing (or fitting) $m_i$, which means that $\varrho _D(v)=m_i(v)$
for every node $v\in V$ where $\varrho _D(v)$ denotes the number of
arcs of $D$ with head $v$.  Often we call a function $m_i$ an {\bf
in-degree specification} or {\bf sequence} or {\bf prescription}.  An
out-degree specification $m_o$ is defined analogously, and a pair
$(m_o,m_i)$ of functions is simply called a {\bf
degree-specification}.

For any function $m:S\rightarrow \mathbb{R}$, the set-function
$\widetilde m$ is defined by $\widetilde m(X):=\sum [m(v):v\in X]$ for
$X\subseteq S$, and we shall use this tilde-notation $\widetilde m$
throughout the paper.  In order to realize $(m_o,m_i)$ with a digraph,
it is necessary to require that $\widetilde m_o(V)=\widetilde m_i(V)$
since both $\widetilde m_o(V)$ and $\widetilde m_i(V)$ enumerate the
total number of arcs of a realizing digraph $D$.  This common value
will be denoted by $\gamma $, that is, our assumption throughout is
that \begin{equation}\widetilde m_o(V)=\widetilde m_i(V)=\gamma . \label{(gamma)}
\end{equation}

The following result was proved by Ore \cite{Ore56} in a slightly
different but equivalent form.

\begin{theorem}[Ore] \label{Ore} A digraph $H=(V,F)$ has a subgraph fitting
$(m_o,m_i)$ if and only if \begin{equation}\widetilde m_i(X) + \widetilde m_o(Z) -
d_H(Z,X) \leq \gamma \ \hbox{for every}\ X,Z\subseteq V, \label{(Ore)}
\end{equation} where $d_H(Z,X)$ denotes the number of arcs $uv\in F$ with $u\in
Z$ and $v\in X$.  \end{theorem}

This immediately implies the following characterization (\cite{Ore56},
see also \cite{Ford-Fulkerson}).

\begin{theorem} [Ore] \label{kinai.alap} Let $(m_o,m_i)$ be a
degree-specification meeting \eqref{(gamma)}.  \medskip

\noindent {\bf (A)} There always exists a digraph realizing
$(m_o,m_i)$.

\medskip \noindent {\bf (B)} There is a loopless digraph realizing
$(m_o,m_i)$ if and only if \begin{equation}\hbox{$m_i(v)+m_o(v)\leq \gamma $ for
every $v\in V$.}\ \label{(kinai.loopless)} \end{equation}

\noindent {\bf (C)} \ There is a simple digraph realizing $(m_o, m_i)$
if and only if \begin{equation}\widetilde m_i(X) + \widetilde m_o(Z) - \vert
X\vert \vert Z\vert + \vert X\cap Z\vert \leq \gamma \ \hbox{for
every}\ X,Z\subseteq V. \label{(kinai.alap)} \end{equation} Moreover, it
suffices to require the inequality in \eqref{(kinai.alap)} only for
its special case when $X$ consists of the $h$ largest values of $m_i$
and $Z$ consists of the $j$ largest values of $m_o$ \ $(h=1,\dots ,n,
\ j=1,\dots ,n)$.  \end{theorem}

Note that \eqref{(kinai.loopless)} follows from \eqref{(kinai.alap)} by
taking $X=\{v\}$ and $Z=\{v\}$.  We also remark that Part (A) can be
proved directly by a simple greedy algorithm:  build up a digraph by
adding arcs $uv$ one by one as long as there are (possibly not
distinct) nodes $u$ and $v$ with $\varrho (v)<m_i(v)$ and $\delta
(v)<m_o(v)$.  Also, Part (B) immediately follows by the following {\bf
loop-reducing technique}.  Let $D$ be a digraph fitting $(m_o,m_i)$
and suppose that there is a loop $e=vv$ sitting at $v$.  Condition
\eqref{(kinai.loopless)} implies that there must be an arc $f=xy$ with
$x\not =v\not =y$ (but $x=y$ allowed).  By replacing $e$ and $f$ with
arcs $xv$ and $vy$, we obtain a digraph fitting $(m_o,m_i)$ that has
fewer loops than $D$ has.  For later purposes, we remark that the
loop-reducing procedure does not decrease the in-degree of any subset
of nodes.

We call a digraph $D=(V,A)$ {\bf strongly connected} or just {\bf
strong} if $\varrho _D(X)\geq 1$ whenever $\emptyset \subset X\subset
V$.  More generally, $D$ is {\bf $k$-edge-connected} if $\varrho
_D(X)\geq k$ whenever $\emptyset \subset X\subset V$.  $D$ is {\bf
$k$-node-connected} or just {\bf $k$-connected} if $k\leq \vert V\vert
-1$ and the removal of any set of less than $k$ nodes leaves a strong
digraph.

One may be interested in characterizing degree-sequences of
$k$-edge-connected and $k$-node-connected digraphs.  We will refer to
this kind of problems as synthesis problems.  The more general
augmentation problem consists of making an initial digraph
$D_0=(V,A_0)$ $k$-edge- or $k$-node-connected by adding a
degree-specified digraph.  Clearly, when $A_0$ is empty we are back at
the synthesis problem.  The augmentation problem was solved for
$k$-edge-connectivity in \cite{FrankJ23} and for $k$-node-connectivity
in \cite{FrankJ31}, but in both cases the augmenting digraphs $D$ were
allowed to have loops or parallel arcs.  The same approach rather
easily extends to the case when $D$ is requested to be loopless but
treating simplicity is significantly more difficult.

The goal of the present paper is to investigate these degree-specified
augmentation and synthesis problems when simplicity is expected.  In
the augmentation problem, this means actually two possible versions
depending on whether the augment{\bf ing} digraph $D$ or else the
augment{\bf ed} digraph $D_0+D$ is requested to be simple.  Clearly,
when $D_0$ has no arcs (the synthesis problem) the two versions
coincide.  We will consider both variations.

An early result of this type is due to Beineke and Harary
\cite{Beineke-Harary} who characterized degree-sequences of loopless
strongly connected digraphs.  In a recent work Hong, Liu, and Lai
\cite{Hong-Liu-Lai} characterized degree-sequences of simple strongly
connected digraphs.  In order to generalize conveniently their result,
we formulate it in a slightly different but equivalent form.

\begin{theorem} [Hong, Liu, and Lai] \label{kinai} Suppose that there is a
simple digraph fitting a degree-specification $(m_o, m_i)$ \ (that is,
\eqref{(kinai.alap)} holds).  There is a strongly connected simple
digraph fitting $(m_o, m_i)$ if and only if \begin{equation}\widetilde m_o(Z) +
\widetilde m_i(X) -\vert X\vert \vert Z\vert + 1 \leq \gamma
\label{(kinai.felt)} \end{equation} holds for every pair of disjoint subsets
$X,Z\subset V$ with $X\cup Z\not =\emptyset $. Moreover, it suffices
to require \eqref{(kinai.felt)} only in the special case when $X$
consists of the $h$ largest values of $m_i$ and $Z$ consists of the
$j$ largest values of $m_o$ $(j=0,1,\dots ,n, \ h=0,1,\dots ,n, \
1\leq j+h\leq n)$.  \end{theorem}

We are going to extend this result in two directions.  In the first
one, degree-specifications are characterized for which there is a
simple realizing digraph $D$ whose addition to an initial
$(k-1)$-edge-connected digraph $D_0$ results in a $k$-edge-connected
digraph $D_0+D$.  The general problem of augmenting an arbitrary
initial digraph $D_0$ with a degree-specified simple digraph to obtain
a $k$-edge-connected digraph remains open even in the special case
when $D_0$ has no arcs at all.  That is, the synthesis problem of
characterizing degree-sequences of simple $k$-edge-connected digraphs
remains open for $k\geq 2$.

Our second generalization of Theorem \ref{kinai} provides a
characterization of degree-sequences of simple $k$-node-connected
digraphs.  We also solve the more general degree-specified
node-connectivity augmentation problem when the augmented digraph is
requested to be simple.  It is a bit surprising that node-connectivity
augmentation problems are typically more complex than their
edge-connectivity counterparts and yet an analogous characterization
for the general $k$-edge-connected case, as indicated above, remains
open.

In the proof of both extensions, we rely on the following general
result of Frank and Jord\'an \cite{FrankJ31}.

\begin{theorem} [Supermodular arc-covering, bi-set function version]
\label{arc-covering} Let $p$ be a positively crossing supermodular
bi-set function which is zero on trivial bi-sets.  The minimum number
of arcs of a loopless digraph covering $p$ is equal to $\max \{\sum
[p(B):  \ B \in {\cal I}]\}$ where the maximum is taken over all
independent families $\cal I$ of bi-sets.  There is an algorithm for
crossing supermodular $p$, which is polynomial in $\vert V\vert $ and
in the maximum value of $p(B)$, to compute the optima.  \end{theorem}

One way to obtain a non-negative positively crossing supermodular
function is taking a crossing supermodular function and increase
its negative values to zero, but not every non-negative positively crossing supermodular function arises in this way.  In deriving
applications, it is simpler to work with the more general notion but
it should be emphasised that no polynomial algorithm is known when $p$
is positively crossing supermodular, and the existing algorithms
work only for crossing supermodular functions.

The algorithm described in \cite{FrankJ31} relies on the ellipsoid
method and on a subroutine to minimize a submodular function given on
a ring-family (see, Schrijver's algorithm in \cite{Schrijver2000}).
V\'egh and Bencz\'ur \cite{Benczur-Vegh} developed a purely
combinatorial algorithm for optimal directed node-connectivity
augmentation problem, an important special case of Theorem
\ref{arc-covering}.  Though not mentioned explicitly in the paper of
V\'egh and Bencz\'ur, their algorithm can be extended without much
effort to the general case described in Theorem \ref{arc-covering}
when $p$ is crossing supermodular.  This extended algorithm
relies on an oracle for minimizing submodular functions.  We should
emphasize that the algorithm of V\'egh and Bencz\'ur is particularly
intricate and it is a natural goal to develop simpler algorithms for
the special cases considered in the present work.

This theorem was earlier used to solve several connectivity
augmentation problems.  It should be, however, emphasized that even
this general framework did not allow to handle simplicity.  Even
worse, there is no hope to extend Theorem \ref{arc-covering} so as to
characterize minimal {\em simple} digraphs covering $p$ since this
problem, by relying on a result by D\"urr et al.  \cite{DGM}, was
shown in \cite{Berczi-Frank16a} (Theorem 12) to include {\bf
NP}-complete special cases.  In \cite{Berczi-Frank16a} and
\cite{Berczi-Frank16b}, we developed other applications of the
supermodular arc-covering theorem when simplicity could be guaranteed.

\medskip A main feature of the present approach to manage the
above-mentioned special cases (when simplicity of the augmenting
digraph is an expectation) is that, though Theorem \ref{arc-covering}
remains a fundamental starting point, relatively tedious additional
work is needed.  (The complications may be explained by the fact that
some special cases are {\bf NP}-complete while others are in {\bf
NP}$\cap $co-{\bf NP}.)

There are actually two issues here to be considered.  The first one is
to develop techniques for embedding special simplicity-requesting
connectivity augmentation problems into the framework of Theorem
\ref{arc-covering}.  When this is successful, one has to resolve a
second difficulty stemming from the somewhat complicated nature of an
independent family of bi-sets in Theorem \ref{arc-covering}.

To demonstrate this second obstacle, consider the following digraph
$D_0=(V,A_0)$ with a particularly simple structure.  Let $e=uv$ be an
arc of $D_0$ if $u\in Z$ or $v\in X$, where $Z$ and $X$ are two
specified disjoint subsets of $V$ with $\vert X\vert =\vert Z\vert
<k$.  An earlier direct consequence of Theorem \ref{arc-covering}
(formulated in Section \ref{node-con} as Theorem \ref{node.augment})
does provide a formula for the minimum number of new arcs whose
addition to any initial digraph results in a $k$-connected digraph.
But to prove that this minimum for our special digraph $D_0$ is
actually equal to the total out-deficiency $\sum [(k-\delta _0(v))\sp
+:  v\in V]$ of the nodes of $D_0$ is rather tricky or tedious.

From an algorithmic point of view, it should be noted that the
original proof of Hong et al.  \cite{Hong-Liu-Lai} gives rise to a
polynomial algorithm which is purely combinatorial.  The present
approach makes use of the supermodular arc-covering theorem.  Since
\cite{FrankJ31} and \cite{Benczur-Vegh} describe polynomial algorithm
to compute the optima in question, there are polynomial algorithms for
finding the simple degree-specified digraphs with the prescribed
connectivity properties.  However, the algorithm in \cite{FrankJ31}
relies on the ellipsoid method while the general algorithm of V\'egh
and Bencz\'ur is particularly complex.  Therefore developing a simple
combinatorial algorithm for our cases requires further investigations.

Also, when, instead of exact degree-specifications, upper and lower
bounds are prescribed for the in-degrees and out-degrees, the problem
of characterizing the existence of simple $k$-node-connected and
degree-constrained digraphs, even in the special case $k=1$, remains a
challenging research task for the future.

\subsection{Notions and notation} \label{Notions}

For a number $x$, let $x\sp += \max\{0,x\}$.  For a function
$m:V\rightarrow \mathbb{R}$ and for $X\subseteq V$, let $\widetilde
m(X):=\sum [m(v):v\in X]$.  For a set-function $p$ and a family $\cal
F$ of sets, let $\widetilde p({\cal F}):  =\sum [p(Z):Z\in {\cal F}].$
For a family $\cal T$ of sets, let $\cup {\cal T}$ denote the union of
the members of $\cal T$.

Two subsets of a ground-set $V$ are said to be {\bf co-disjoint} if
their complements are disjoint.  By a {\bf partition} of a ground-set
$V$, we mean a family of disjoint subsets of $V$ whose union is $V$.
A {\bf subpartition} of $V$ is a partition of a subset of $V$.  A {\bf
co-partition} (resp., {\bf co-subpartition}) of $V$ is a family of
subsets arising from a partition (subpartition) of $V$ by
complementing each of its member.  For a subpartition ${\cal
T}=\{T_1,\dots ,T_q\}$, we always assume that its members $T_i$ are
non-empty but $\cal T$ is allowed to be empty (that is, $q=0$).

Two subsets $X,Y\subseteq V$ are {\bf crossing} if none of $X-Y,Y-X,
X\cap Y, V-(X\cup Y)$ is empty.  A family of subsets is {\bf
cross-free} if it contains no two crossing members.  A family $\cal F$
of subsets is {\bf crossing} if both $X\cap Y$ and $X\cup Y$ belong to
$\cal F$ whenever $X$ and $Y$ are crossing members of $\cal F$.

When $X$ is a subset of $V$, we write $X\subseteq V$, while $X\subset
V$ means that $X$ is a proper subset.  The standard notation $A
\setminus B$ for set difference will be replaced by $A-B$.  When it
does not cause any confusion, we do not distinguish a one-element set
$\{v\}$ (often called a {\bf singleton}) from its only element $v$,
and we use the notation $v$ for the singleton as well.  For example,
we write $V-v$ rather than $V-\{v\}$, and $V+v$ stands for $V\cup
\{v\}$.  In some situations, however, the formally precise $\{v\}$
notation must be used.  For example, an arc $e=uv$ in a digraph is
said to enter a node $v$ even if $e$ is a loop (that is, $u=v$) while
$e$ enters the singleton $\{v\}$ only if $u\in V-v$.  That is a loop
sitting at $v$ enters $v$ but does not enter $\{v\}$.  Therefore the
in-degree $\varrho (v)$ (the number of arcs entering $v$) is equal to
the in-degree $\varrho (\{v\})$ plus the number of loops sitting at
$v$.

In a digraph $D=(V,A)$, an arc $uv$ {\bf enters} a subset $X\subseteq
V$ or {\bf leaves} $V-X$ if $u\in X, \ v\in V-X$.  The in-degree
$\varrho _D(X)=\varrho _A(X)$ of a subset $X\subseteq V$ is the number
of arcs entering $X$ while the out-degree $\delta _D(X)=\delta _A(X)$
is the number of arcs leaving $X$.  Two arcs of a digraph are {\bf
parallel} if their heads coincide and their tails coincide.  The
oppositely directed arcs $uv$ and $vu$ are not parallel.  We call a
digraph {\bf simple} if it has neither loops nor parallel arcs.

By the {\bf complete digraph} $D\sp *=(V,A\sp *)$, we mean the simple
digraph on $V$ in which there is one arc from $u$ to $v$ for each
ordered pair $(u,v)$ of distinct nodes, that is, $D\sp *$ has $n(n-1)$
arcs.  For two subsets $Z,X\subseteq V$, let $D\sp *[Z,X]$ denotes the
subgraph of $D\sp *$ consisting of those arcs $uv$ for which $u\in Z$
or $v\in X$.  Then $D\sp *[Z,X]$ has $\vert Z\vert (n-1)+(n-\vert
Z\vert )\vert X\vert -\vert X-Z\vert $ arcs.

For two digraphs $D_0=(V,A_0)$ and $D=(V,A)$ on the same node-set,
$D_0+D=(V,A_0+A)$ denotes the digraph consisting of the arcs of $D_0$
and $D$.  That is, $D_0+D$ has $\vert A_0\vert +\vert A\vert $ arcs.

A digraph $D$ {\bf covers} a family $\cal K$ of subsets if $\varrho
_D(K)\geq 1$ for every $K\in {\cal K}$.  A digraph $D$ {\bf covers} a
set-function $p$ on $V$ if $\varrho _D(K)\geq p(K)$ for every
$K\subseteq V$.

By a {\bf bi-set} we mean a pair $B=(B_O,B_I)$ of subsets for which
$B_I\subseteq B_O$.  Here $B_O$ and $B_I$ are the outer set and the
inner set of $B$, respectively.  A bi-set is {\bf trivial} if
$B_I=\emptyset $ or $B_O=V$.  The set $W(B):=B_O-B_I$ is the {\bf
wall} of $B$, while $w(B)=\vert W(B)\vert $ is its {\bf wall-size}.
Two bi-sets $B$ and $C$ are {\bf comparable} if $B_I\subseteq C_I$ and
$B_O\subseteq C_O$ or else $B_I\supseteq C_I$ and $B_O\supseteq C_O$.
The {\bf meet} of two bi-sets $B$ and $C$ is defined by $B\sqcap
C=(B_O\cap C_O,B_I\cap C_I)$ while their {\bf join} is $B\sqcup
C=(B_O\cup C_O,B_I\cup C_I)$.  Note that $w(B)$ is a modular function
in the sense that $w(B)+w(C) = w(B\sqcap C) + w(B\sqcup C)$.  Two
bi-sets $B$ and $C$ are {\bf crossing} if $B_O\cup C_O\not =V$,
$B_I\cap C_I\not =\emptyset $, and they are not comparable.  A family
$\cal B$ of bi-sets is {\bf crossing} if both $B\sqcap C$ and $B\sqcup
C$ belong to $\cal B$ whenever $B$ and $C$ are crossing members of
$\cal B$.  A bi-set function $p$ is {\bf positively crossing
supermodular} if $$p(B)+p(C)\leq p(B\sqcap C) + p(B\sqcup C)$$
whenever $p(B)>0, p(C)>0$, $B$ and $C$ are crossing bi-sets.

An arc $e$ {\bf enters} (or {\bf covers}) a bi-set $B$ if $e$ enters
both $B_O$ and $B_I$.  The {\bf in-degree} $\varrho (B)$ of a bi-set
$B$ is the number of arcs entering $B$.  Two bi-sets are {\bf
independent} if no arc can cover both, which is equivalent to
requiring that their inner sets are disjoint or their outer sets are
co-disjoint.  A family of bi-sets is {\bf independent} if their
members are pairwise independent.  Given a digraph $D=(V,A)$, a bi-set
$B=(B_O,B_I)$ is {\bf $D$-one-way} or just {\bf one-way} if no arc of
$D$ covers $B$.

\medskip %{\tiny {\bf directory:  fenyo, file:  novel1 \today}}

\section{Edge-connectivity} \label{edge-con}

The degree-specified augmentation problem for $k$-edge-connectivity
was shown by the second author \cite{FrankJ23} to be equivalent to
Mader's directed splitting off theorem \cite{Mader82}.

\begin{theorem} [\cite{FrankJ23}] \label{edge.augment.m.loopless} An initial
digraph $D_0=(V,A_0)$ can be made $k$-edge-connected by adding a
digraph $D=(V,A)$ fitting $(m_o,m_i)$ if and only if $\widetilde
m_i(X)+\varrho _{D_0}(X)\geq k$ and $\widetilde m_o(X)+\delta
_{D_0}(X)\geq k$ hold for every subset $\emptyset \subset X\subset V$.
If, in addition, \eqref{(kinai.loopless)} holds, then $D$ can be
chosen loopless.  \end{theorem}

The second part immediately follows from the first one by applying the
loop-reducing technique mentioned in Section \ref{intro} since
loop-reduction never decreases the in-degree of a subset.  Though the
problem when simplicity of $D$ is requested remains open even in the
special case when $D_0$ has no arcs, we are able to prove the
following straight extension of Theorem \ref{kinai}.

\begin{theorem} \label{erosnov} Suppose that there is a simple digraph
fitting a degree-specification $(m_o, m_i)$ \ (that is,
\eqref{(kinai.alap)} holds).  A digraph $D_0$ can be made strongly
connected by adding a simple digraph fitting a degree-specification
$(m_o,m_i)$ if and only if inequality \eqref{(kinai.felt)} holds for
every pair of disjoint subsets $X,Z\subset V$ for which there is a
non-empty, proper subset $K$ of $V$ so that $\varrho _{D_0}(K)=0$ and
$Z\subseteq K\subseteq V-X$.\end{theorem}

This theorem is just a special case of the following.

\begin{theorem} \label{elofnov} Suppose that there is a simple digraph
fitting a degree-specification $(m_o, m_i)$.  A $(k-1)$-edge-connected
digraph $D_0$ can be made $k$-edge-connected by adding a simple
digraph fitting a degree-specification $(m_o,m_i)$ if and only if
inequality \eqref{(kinai.felt)} holds for every pair of disjoint
subsets $X,Z\subset V$ for which there is a non-empty, proper subset
$K$ of $V$ so that $\varrho _{D_0}(K)=k-1$ and $Z\subseteq K\subseteq
V-X$.\end{theorem}

Since the family of subsets of in-degree $k-1$ in a
$(k-1)$-edge-connected digraph $D_0$ is a crossing family, the
following result immediately implies Theorem \ref{elofnov}.

\begin{theorem} \label{kinai.main} Let $\cal K$ be a crossing family of
non-empty proper subsets of $V$.  Suppose that there is a simple
digraph fitting a degree specification $(m_o, m_i)$, that is,
\eqref{(kinai.alap)} holds.  There is a simple digraph fitting
$(m_o,m_i)$ which covers $\cal K$ if and only if \begin{equation}\widetilde m_o(Z)
+ \widetilde m_i(X) -\vert X\vert \vert Z\vert +1\leq \gamma
\label{(kinai.ujfelt)} \end{equation} holds for every pair of disjoint subsets
$X,Z\subset V$ for which there is a member $K\in {\cal K}$ with
$Z\subseteq K\subseteq V-X$.  \end{theorem}

\proof{Proof.}
Suppose that there is a requested digraph $D$.  By the
simplicity of $D$, there are at most $\vert X\vert \vert Z\vert $ arcs
from $Z$ to $X$.  Therefore the total number of arcs with tail in $Z$
or with head in $X$ is at least $\widetilde m_o(Z) + \widetilde m_i(X)
-\vert X\vert \vert Z\vert $. Moreover at least one arc enters $K$ and
such an arc neither leaves an element of $Z$ nor enters an element of
$X$, from which we obtain that $\widetilde m_o(Z) + \widetilde m_i(X)
-\vert X\vert \vert Z\vert +1\leq \gamma $, that is,
\eqref{(kinai.ujfelt)} is necessary.

To prove sufficiency, observe first that the theorem is trivial if
$n:=\vert V\vert \leq 2$ so we assume that $n\geq 3$.  We need some
further observations.

\begin{claim} \label{kinai.claim1} \begin{equation}\hbox{ $\widetilde m_i(K)\geq 1$ and
$\widetilde m_o(V-K)\geq 1$ holds for every $K\in {\cal K}$.}\ \end{equation} In
particular, if $\{v\}\in {\cal K}$ for some $v\in V$, then $m_i(v)\geq
1$, and if $V-u\in {\cal K}$ for some $u\in U$, then $m_o(u)\geq 1$.
\end{claim}

\proof{Proof.}
$\widetilde m_i(K)\geq 1$ follows by applying
\eqref{(kinai.ujfelt)} to $Z=\emptyset $ and $X=V-K$, while $\widetilde m_o(V-K)\geq 1$ follows with the choice $X=\emptyset $ and $Z=K$.
$\bullet$\endproof

\medskip This claim immediately implies the following.

\begin{claim} \label{kinai.gamma} $\cal K$ has at most $\gamma $ pairwise
disjoint and at most $\gamma $ pairwise co-disjoint members.\end{claim}

\begin{claim} \label{mov-miv} $m_o(v)\leq n-1$ and $m_i(v)\leq n-1$ for every
$v\in V.$ \end{claim}

\proof{Proof.}
By applying \eqref{(kinai.alap)} to $X=\{v\}$ and $Z=V$, one
gets $m_i(v) + \widetilde m_o(V) - 1\cdot\vert V\vert + \vert
\{v\}\vert \leq \gamma $, that is, $m_i(v)\leq n-1$, and $m_o(v)\leq
n-1$ is obtained analogously by choosing $X=V$ and $Z=\{v\}$.
$\bullet$\endproof \medskip

Define a bi-set function $p$ as follows.  Let $p(Y_O,Y_I)$ be zero
everywhere apart from the next three types of bi-sets.  \medskip

\noindent Type 1:  \ For $K\in {\cal K}$ with $1<\vert K\vert <n-1$,
let $p(K,K)=1$.

\medskip \noindent Type 2:  \ For $u\in V$, let $p(V-u,V-u)=m_o(u)$.

\medskip

\noindent Type 3:  \ For $v\in Y\subset V$, let $p(Y,\{v\})= m_i(v) -
(\vert Y\vert -1).$

\medskip

Note that the role of $m_o$ and $m_i$ is not symmetric in the
definition of $p$.  Since $n\geq 3$, each bi-set $B$ with positive
$p(B)$ belongs to exactly one of the three types.

\begin{claim} \label{kinai.crossing} The bi-set function $p$ defined above is
positively crossing supermodular.  \end{claim}

\proof{Proof.}Let $B=(B_O,B_I)$ and $C=(C_O,C_I)$ be two crossing bi-sets
with $p(B)>0$ and $p(C)>0$.  Then neither of $B$ and $C$ is of Type 2.
Suppose first that both $B$ and $C$ are of Type 1. Observe that if
$K=\{v\}\in \cal K$ for some $v\in V$, then $(K,K)$ is of Type 3 and
hence $p(K,K)= m_i(v)\geq 1$ by Claim \ref{kinai.claim1}.  Similarly,
if $K=V-u\in \cal K$ for some $u\in V$, then $(K,K)$ is of Type 2 and
hence $p(K,K)= m_o(u)\geq 1$.  Therefore the supermodular inequality
in this case follows from the assumption that the set-system $\cal K$
is crossing.

If both $B$ and $C$ are of Type 3, then $B_I=\{v\}=C_I$ for some $v\in
V$, and in this case the supermodular inequality holds actually with
equality.

Finally, let $B$ be of Type 1 and let $C$ be of Type 3. Then
$B_I=K=B_O$ for some $K\in {\cal K}$ with $1<\vert K\vert <n-1$ and
$C_I=\{v\}$ for some $v\in K$.  Observe that $B\sqcup C$ does not
belong to any of the three types and hence $p(B\sqcup C)=0$.
Furthermore, since $B\sqcap C$ is of Type 3, we have $p(B\sqcap
C)=m_i(v) -(\vert B_O\cap C_O\vert -1)$.

Since $(K,K)$ and $(C_O,\{v\})$ are not comparable, $\vert C_O\cap
K\vert \leq \vert C_O\vert -1$ and therefore $p(B) + p(C) = 1 +
[m_i(v)- (\vert C_O\vert -1)] \leq [m_i(v) -(\vert B_O\cap C_O\vert
-1)] + 0 = p(B\sqcap C) + p(B\sqcup C).$ $\bullet$\endproof \medskip

It follows from the definition of $p$ that every digraph covering $p$
must have at least $\gamma $ arcs.  \medskip

\noindent {\bf Case 1}.  \ There is a loopless digraph $D=(V,A)$ with
$\gamma $ arcs covering $p$.  \medskip

Now $\gamma =\vert A\vert =\sum [\varrho _A(v):v\in V]\geq \sum
[m_i(v):v\in V] = \gamma $ from which $\varrho _A(v)=m_i(v)$ follows
for every $v\in V$.  Analogously, we get $\delta _A(v)= \varrho
_A(V-v)= m_o(v)$ for every $v\in V$.  By the definition of $p$, it
also follows that $D$ covers $\cal K$.

\begin{claim} \label{kinai.simple} $D$ is simple.  \end{claim}

\proof{Proof.}Suppose, indirectly, that $D$ has two parallel arcs $e$ and $f$
from $u$ to $v$.  Consider the bi-set $(Y,\{v\})$ for $Y=\{u,v\}$.  We
have $p(Y,\{v\}) \leq \varrho (v)-2 =m_i(v)-2 = p(Y,\{v\})-1$, a
contradiction.  $\bullet$\endproof

\medskip We can conclude that in Case 1 the digraph requested by the
theorem is indeed available.  \medskip

\noindent {\bf Case 2}.  \ The minimum number of arcs of a loopless
digraph covering $p$ is larger than $\gamma $. \medskip

Theorem \ref{arc-covering} implies that in Case 2 there is an
independent family ${\cal I}$ of bi-sets for which $\widetilde p({\cal
I})>\gamma $. Then $\cal I$ partitions into three parts according to
the three possible types its members belong to.  Therefore we have a
subset ${\cal F}\subseteq {\cal K}$, a subset $Z\subseteq V$, and a
family ${\cal B}=\{(Y,\{v\}):  v\in Y\subset V\}$ of bi-sets so that
$$ {\cal I} = \{(K,K):  K\in {\cal F}\} \cup \{(V-z,V-z):  z\in Z\}
\cup {\cal B}, \hbox{and}\ $$ \begin{equation}\vert {\cal F}\vert + \widetilde
m_o(Z) + \sum [m_i(v) - (\vert Y\vert -1) :  (Y,\{v\})\in {\cal B}] =
\widetilde p({\cal I}) >\gamma . \label{(kinai.nagyobb)} \end{equation}

We may assume that $\vert \cal B\vert $ is as small as possible.

\begin{claim} \label{kinai.egybeta} There are no two members $(Y,\{v\})$ and
$(Y',\{v\})$ of $\cal B$ with the same inner set $\{v\}$.\end{claim}

\proof{Proof.}If indirectly there are two such members, then $Y\cup Y'=V$ by
the independence of $\cal I$.  If we replace the two members
$(Y,\{v\})$ and $(Y',\{v\})$ of $\cal I$ by the single bi-set $(Y\cap
Y',\{v\})$, then the resulting family ${\cal I}'$ is also independent
since any arc covering $(Y\cap Y',\{v\})$ covers at least one of
$(Y,\{v\})$ and $(Y',\{v\})$.  Furthermore, \begin{equation}p(Y,\{v\}) +
p(Y',\{v\}) = m_i(v) - (\vert Y\vert -1) + m_i(v) - (\vert Y'\vert -1)
\label{(kinai.becsles1)} \end{equation} and \begin{equation}p(Y\cap Y',\{v\}) = m_i(v) -
(\vert Y\cap Y'\vert -1).  \label{(kinai.becsles2)} \end{equation}

We claim that $p(Y,\{v\}) + p(Y',\{v\}) \leq p(Y\cap Y',\{v\})$.
Indeed, this is equivalent to $m_i(v)\leq (\vert Y\vert -1) + (\vert
Y'\vert -1) - (\vert Y\cap Y'\vert -1)$, that is, $m_i(v)\leq \vert
V\vert -1$, but this holds by Claim \ref{mov-miv}.  Consequently,
$\widetilde p({\cal I}')\geq \widetilde p({\cal I})>\gamma $,
contradicting the minimality of $\vert \cal B\vert $. $\bullet$\endproof
\medskip

Let $X:=\{v:$ there is a bi-set $(Y,\{v\})\in {\cal B}\}$.  For any
element $v\in X$, Claim \ref{kinai.egybeta} implies that the (outer)
set $Y$ for which $(Y,\{v\})\in {\cal B}$ is uniquely determined, and
it will be denoted by $Y_v$.  Now \eqref{(kinai.nagyobb)} transforms to
\begin{equation}\vert {\cal F}\vert + \widetilde m_o(Z) + \sum [m_i(v) - (\vert
Y_v\vert -1) :  v\in X ] = \widetilde p({\cal I}) >\gamma .
\label{(kinai.nagyobb.2)} \end{equation}

We may assume that $\widetilde p({\cal I})$ is as large as possible,
and modulo this, $\vert {\cal F}\vert $ is minimal.

\begin{claim} \label{kinai.Kv} For $v\in X-Z$, one has $Y_v=Z+v$, and for
$v\in X\cap Z$, one has $Y_v=Z$.  \end{claim}

\proof{Proof.}The independence of $\cal I$ implies $Z\subseteq Y_v$.
Suppose, indirectly, that there is an element $u\in (Y_v-v)-Z$.
Replace the member $(Y_v,\{v\})$ of $\cal I$ by $(Y_v-u,\{v\})$.
Since $p(Y_v-u,\{v\})=p(Y_v)+1,$ the resulting system ${\cal I}'$ is
not independent, therefore there exists a member $(Y_O,Y_I)$ of $\cal
I$ that is covered by arc $uv$.  Since $\cal I$ is independent, this
member is unique.  As $u\not \in Z$, Claim \ref{kinai.egybeta} implies
that $(Y_O,Y_I)$ must be in $\cal F$, that is, $(Y_O,Y_I)=(K,K)$ for
some $K\in {\cal K}$.  By leaving out $(K,K)$ from ${\cal I}'$, we
obtain an independent ${\cal I}''$ for which $\widetilde p({\cal
I}'')=\widetilde p({\cal I})$, contradicting the minimal choice of
$\cal F$.  $\bullet$\endproof

\medskip Due to these claims, Condition \eqref{(kinai.nagyobb)} reduces
to $\vert {\cal F}\vert + \widetilde m_o(Z) + \widetilde m_i(X) -
[\vert X\cap Z\vert (\vert Z\vert -1)+ \vert X-Z\vert \vert Z\vert ]
>\gamma $, which is equivalent to \begin{equation}\vert {\cal F}\vert + \widetilde
m_o(Z) + \widetilde m_i(X) - \vert X\vert \vert Z\vert + \vert X\cap
Z\vert >\gamma . \label{(kinai.nagyobb.3)} \end{equation}

\begin{claim} \ $X\cup Z\not =\emptyset $. \end{claim}

\proof{Proof.}If $X\cup Z=\emptyset $, then \eqref{(kinai.nagyobb.3)} reduces
to $\vert \cal F\vert >\gamma $. It is an easy observation that the
members of the independent $\cal F$ are either pairwise disjoint or
pairwise co-disjoint, contradicting Claim \ref{kinai.gamma}.  $\bullet
$

\begin{claim} \label{Znemures} $Z\not =\emptyset $. \end{claim}

\proof{Proof.}If indirectly $Z=\emptyset $, then $X\not =\emptyset $ in which
case \eqref{(kinai.nagyobb.3)} reduces to

\begin{equation}\vert {\cal F}\vert + \widetilde m_i(X) >\gamma .
\label{(Znemures)} \end{equation}

By Claim \ref{kinai.Kv}, we have $Y_v=\{v\}$ for every element $v\in
X$.  The independence of $\cal I$ implies that all the members of
$\cal F$ are disjoint from $X$.  Hence ${\cal F}=\{P_1,\dots ,P_q\}$
is a subpartition of $V-X$.  Now \eqref{(Znemures)} is equivalent to
$q>\widetilde m_i(V-X)$, and by Claim \ref{kinai.claim1} and $m_i\geq
0$, we have $$q>\widetilde m_i(V-X) \geq \widetilde m_i(\cup {\cal F})
= \sum _{j=1}\sp q \widetilde m_i(P_j) \geq q,$$ a contradiction.
$\bullet$\endproof

\begin{claim} \label{Xnemures} $X\not =\emptyset $. \end{claim}

\proof{Proof.}
Suppose, indirectly, that $X=\emptyset .$ Then
\eqref{(kinai.nagyobb.3)} reduces to \begin{equation}\vert {\cal F}\vert +
\widetilde m_o(Z) >\gamma . \label{(kinai.Z)}\end{equation}

The independence of $\cal I$ implies that all the members of ${\cal
F}$ include $Z$ and these members are pairwise co-disjoint.  Let
${\cal P}=\{P_1,\dots ,P_q\}$ consist of the complements of the
members of $\cal F$ (that is, ${\cal F}=\{V-P_1,\dots ,V-P_q\}$).
Then $\cal P$ is a subpartition of $V-Z$.

By \eqref{(kinai.Z)}, by $m_o\geq 0$, and by Claim \ref{kinai.claim1},
we have $$q=\vert {\cal F}\vert > \gamma - \widetilde m_o(Z) =
\widetilde m_o(V) - \widetilde m_o(Z) \geq \sum _{j=1}\sp q \widetilde
m_o(P_j) \geq q,$$ a contradiction.  $\bullet$\endproof \medskip

We have concluded that $X\not =\emptyset $ and $Z\not =\emptyset $.
\medskip

\noindent {\bf Case A} \ $X\cap Z=\emptyset $. \medskip

On one hand, \eqref{(kinai.alap)} reduces in this case to $\widetilde
m_i(X) +\widetilde m_o(Z) - \vert X\vert \vert Z\vert \leq \gamma .$
On the other hand, \eqref{(kinai.nagyobb.3)} reduces to

\begin{equation}\vert {\cal F}\vert + \widetilde m_o(Z) + \widetilde m_i(X) -
\vert X\vert \vert Z\vert >\gamma . \label{(kinai.nagyobb.4)} \end{equation}

Therefore we must have ${\vert \cal F\vert }\geq 1$.  For $K\in {\cal
F}$, the independence of $\cal I$ implies that $Z\subseteq K$ and
$X\subseteq V-K$. As $X\neq\emptyset $ and $Z\neq\emptyset$, the independence of $\cal I$ also implies that
$\cal F$ cannot have more than one member, that is, ${\vert \cal
F\vert }=1$, and hence $K, X, Z$ violate \eqref{(kinai.ujfelt)}.

\medskip

\noindent {\bf Case B} \ $X\cap Z\not =\emptyset $. \medskip

If $K\in {\cal F}$, then $Z\subseteq K$.  Moreover, for an element
$v\in X\cap Z$, we have $(Y_v,v)=(Z,v)$ by Claim \ref{kinai.Kv}.  The
independence of $(K,K)$ and $(Z,\{v\})$ implies that $Z\cup K=V$, that
is, $K=V$, which is not possible.  Hence ${\cal F}=\emptyset $ and
\eqref{(kinai.nagyobb.3)} reduces to $ \widetilde m_o(Z) + \widetilde
m_i(X) - \vert X\vert \vert Z\vert + \vert X\cap Z\vert >\gamma $
contradicting \eqref{(kinai.alap)}.  This contradiction shows that Case
2 cannot occur, completing the proof of the theorem.  $\bullet$
$\bullet$\endproof

\medskip %{\tiny {\bf directory:  fenyo, file:  novel2, \today}}

\section{Node-connectivity} \label{node-con}

Let $H=(V,F)$ be a simple digraph on $n\geq k+1$ nodes.  Recall that a
bi-set $B=(B_O,B_I)$ was called {\bf $H$-one-way} or just {\bf
one-way} if $\varrho _H(B)=0$, that is, if no arc of $H$ enters both
$B_I$ and $B_O$.  Recall the notation $w(B):=\vert B_O-B_I\vert $. The
following lemma occured in \cite{FrankJ31}.  For completeness, we
include its proof.

\begin{lemma} \label{k-con} The following are equivalent.

\noindent {\em (A1)} \ $H=(V,F)$ is $k$-connected.

\noindent {\em (A2)} \ $\varrho _H(B)+ w(B) \geq k$ for every
non-trivial bi-set $B$.

\noindent {\em (A3)} \ $w(B)\geq k$ for every non-trivial one-way
bi-set $B$.

\noindent {\em (B)} \ There are $k$ openly disjoint $st$-paths in $H$
for every ordered pair of nodes $s,t$.\end{lemma}

\proof{Proof.}(A1) implies (B) by the directed node-version of Menger's
theorem.

(B)$\Rightarrow $(A2).  Let $B=(B_O,B_I)$ be a non-trivial bi-set, let
$s\in V-B_O$ and $t\in B_I.$ By (B), there are $k$ openly disjoint
$st$-paths.  Each of them uses an arc entering $B$ or an element of
the wall $W(B)$ of $B$, from which $\varrho _H(B)+ w(B) \geq k$, that
is, (A2) holds.

(A2)$\Rightarrow $(A3).  Indeed, (A3) is just a special case of (A2).

(A3)$\Rightarrow $(A1).  If (A1) fails to hold, then there is a subset
$Z$ of less than $k$ nodes so that $H-Z$ is not strongly connected.
Let $B_I$ be a non-empty proper subset of $V-Z$ with no entering arc
of $H-Z$ and let $B_O:=Z\cup B_I$.  Then $B=(B_O,B_I)$ is a
non-trivial bi-set with $W(B)=Z$ for which $\varrho _H(B)=0$ and
$w(B)=\vert Z\vert <k$, that is, (A3) also fails to hold.  $\bullet$\endproof

\subsection{Connectivity augmentation:  known results}

Let $D_0=(V,A_0)$ be a starting digraph on $n\geq k+1$ nodes.  The
in-degree and out-degree functions of $D_0$ will be abbreviated by
$\varrho _0$ and $\delta _0$, respectively.  In the connectivity
augmentation problem we want to make $D_0$ $k$-connected by adding new
arcs.  Since parallel arcs and loops do not play any role in
node-connectivity, we may and shall assume that $D_0$ is simple.

In one version of the connectivity augmentation problem, one strives
for minimizing the number of arcs to be added.  In this case, the
optimal augmenting digraph is automatically simple.  The following
result is a direct consequence of Theorem \ref{arc-covering}.

\begin{theorem} [Frank and Jord\'an, \cite{FrankJ31}] \label{node.augment} A
digraph $D_0$ can be made $k$-connected by adding a simple digraph
with at most $\gamma $ arcs if and only if \begin{equation}\sum [k-w(B):  B\in
{\cal I}]\leq \gamma \label{(Frank-Jordan.felt1)} \end{equation} holds for every
independent family $\cal I$ of non-trivial $D_0$-one-way bi-sets.  \end{theorem}

In what follows, ${\cal B}_0$ denotes the family of non-trivial
$D_0$-one-way bi-sets.  For any bi-set $B$, let $p_1(B):=k-w(B)$.
With these terms, \eqref{(Frank-Jordan.felt1)} requires that
$\widetilde p_1({\cal I})\leq \gamma $ for every independent ${\cal
I}\subseteq {\cal B}_0$.

In a related version of the connectivity augmentation problem, the
goal is to find an augmenting digraph $D$ fitting a
degree-specification $(m_o,m_i)$ (meeting \eqref{(gamma)}) so that the
augmented digraph $D_0\sp +:=D_0+D$ is $k$-connected.  The paper
\cite{FrankJ31} described a characterization for the existence of such
a $D$, but in this case the augmenting digraph is not necessarily
simple.  This characterization can also be derived from Theorem
\ref{arc-covering}.

\begin{theorem} \label{node.augment.m} There exists a digraph $D$ fitting
$(m_o,m_i)$ such that $D_0+D$ is $k$-connected if and only if \begin{equation}
\hbox{ $\varrho _0(v)+m_i(v) \geq k$ and $\delta _0(v)+m_o(v)\geq k$
for each $v\in V,$ }\ \label{(Frank-Jordan.pont)} \end{equation} \begin{equation}\widetilde
m_i(Z) \geq \sum [k-w(B):  B\in {\cal I}] \label{(Frank-Jordan.mi)}
\end{equation} holds for every $Z\subseteq V$ and for every independent family
$\cal I$ of non-trivial $D_0$-one-way bi-sets $B$ with $B_I\subseteq
Z$, and \begin{equation}\widetilde m_o(Z) \geq \sum \left [k-w(B):  B\in {\cal I}
\right ] \label{(Frank-Jordan.mo)} \end{equation} holds for every $Z\subseteq V$
and for every independent family $\cal I$ of non-trivial one-way
bi-sets $B$ with $B_O\cup Z=V$.  \end{theorem}

Note that in this theorem both parallel arcs and loops are allowed in
the augmenting digraph $D$.  By using Theorem \ref{arc-covering} and
some standard steps, one can derive the following variation when loops
are excluded.

\begin{theorem} \label{node.augment.m.loopless} There exists a loopless
digraph $D$ fitting $(m_o,m_i)$ such that $D_0+D$ is $k$-connected if
and only if each of \eqref{(kinai.loopless)},
\eqref{(Frank-Jordan.pont)}, \eqref{(Frank-Jordan.mi)}, and
\eqref{(Frank-Jordan.mo)} hold.\end{theorem}

\begin{corollary} If there is a loopless digraph fitting $(m_o,m_i)$ and if
there is a digraph $D$ fitting $(m_o,m_i)$ for which $D_0+D$ is
$k$-connected, then $D$ can be chosen loopless.  \end{corollary}

Note that an analogous statement for $k$-edge-connectivity in Theorem
\ref{edge.augment.m.loopless} follows immediately by applying the
loop-reducing technique, but this approach does not seem to work here
since a loop-reducing step may destroy $k$-node-connectivity.  To see
this, let $D_0$ be the digraph on node-set $V=\{x,y,z,v\}$ with no arc.
Let $k=2$ and let $m_o(x)=2=m_i(x)$, $m_o(y)=2=m_i(y)$, $m_o(z)=2=m_i(z)$, and
$m_o(v)=3=m_i(v)$.  Consider the digraph $D$ with arc set
$\{xy,yx,yz,zy,zv,vz,vx,xv,vv\}$. This digraph is
2-node-connected and fits the degree-specification, but the
loop-reducing technique replaces, say, the two arcs $xy$ and $vv$ by $xv$
and $vy$, and the digraph $D'$ arising in this way is not 2-connected
since $D'-v$ is not strongly connected. Analogously, the same happens when the two arcs $yx$ and $vv$ are replaced by $yv$ and $vx$. Note that the digraph on $V$
with arc set $\{xv,vx,yv,vy,zv,vz,xy,yz,zx\}$ is 2-connected,
loopless, and fits the degree-specification.

\subsection{Degree-specified connectivity augmentation preserving
simplicity}

Our present goal is to solve the degree-specified node-connectivity
augmentation problem when simplicity of the augmented digraph $D_0+D$
is requested.  With the help of a similar approach another natural
variant, when only the simplicity of the augmenting digraph $D$ is
requested, can also be managed.  Theorem \ref{elofnov} provided a
complete answer to this latter problem in the special case $k=1$.

Let $\ol D_0=(V,\ol A_0)$ denote the complementary digraph of $D_0$
arising from the complete digraph $D\sp *=(V,A\sp *)$ by removing
$A_0$, that is, $\ol A_0:=A\sp *-A_0$.  In these terms, our goal is to
find a degree-specified subgraph $D$ of $\ol D_0$ for which $D_0+D$ is
$k$-connected.  Note that in the case when an arbitrary digraph $H$
for possible new arcs is specified instead of $\ol D_0$, the problem
becomes NP-complete even in the special case $k=1$ and $A_0=\emptyset
$ since, for the degree specification $m_o\equiv 1\equiv m_i$, it is
equivalent to finding a Hamiltonian circuit of $H$.

We will show that the problem can be embedded into the framework of
Theorem \ref{arc-covering} in such a way that the augmented digraph
$D_0+D$ provided by this theorem will automatically be simple.  In
this way, we shall obtain a good characterization for the general case
when the initial digraph $D_0$ is arbitrary.  This characterization,
however, will include independent families of bi-sets, and in this
sense it is more complicated than the one given in Theorem \ref{kinai}
for the special case $k=1$.

In the special case when $D_0$ has no arcs at all, that is, when the
goal is to find a simple degree-specified $k$-connected digraph, the general
characterization will be significantly simplified in such a way that
the use of independent bi-set families is completely avoided, and we
shall arrive at a characterization which is a straight extension of
the one in Theorem \ref{kinai} concerning the special case $k=1$.
Recall that the simple digraph $D\sp *[Z,X]$ with node-set $V$ was
defined in Section \ref{Notions} so that $uv$ was an arc for distinct
$u,v\in V$ precisely if $u\in Z$ or $v\in X$.  Also, ${\cal B}_0$ was
introduced above to denote the set of non-trivial $D_0$-one-way
bi-sets.

\begin{theorem} \label{novel3.main} Let $D_0=(V,A_0)$ be a simple digraph
with in- and out-degree functions $\varrho _0$ and $\delta _0$,
respectively.  There is a digraph $D=(V,A)$ fitting $(m_o,m_i)$ for
which $D_0\sp +:=D_0+D$ is simple and $k$-connected if and only if \begin{equation}
\widetilde p_1({\cal F}) + \widetilde m_o(Z) + \widetilde m_i(X)-
d_{\ol A_0}(Z,X) \leq \gamma \label{(felt.F0)} \end{equation} holds for every
pair of subsets $X,Z\subseteq V$ and for every independent family
$\cal F$ of non-trivial bi-sets which are one-way with respect to
$D_0+D\sp *[Z,X]$, where $p_1(B_O,B_I)=k-w(B)$ for $B\in {\cal F}$ and
$d_{\ol A_0}(Z,X)$ denotes the number of arcs $a=zx\in \ol A_0$ for
which $z\in Z$ and $x\in X$.  \end{theorem}

\proof{Proof.}Note that the requirement for the members $B=(B_O,B_I)$ of
${\cal F}$ to be one-way with respect to $D_0+D\sp *[Z,X]$ is
equivalent to requiring that $B\in {\cal B}_0$ and both $Z\subseteq
B_O$ and $X\cap B_I=\emptyset $ hold.

For proving necessity, assume that $D$ is a digraph meeting the
requirements of the theorem.  To see \eqref{(felt.F0)}, observe that
there are $\widetilde m_o(Z)$ arcs of $D$ with tail in $Z$ and there
are $\widetilde m_i(X)$ arcs with head in $X$.  Since $D$ can have at
most $d_{\ol A_0}(Z,X)$ arcs with tail in $Z$ and head in $X$, the
number of arcs of $D$ with tail in $Z$ or head in $X$ is at least
$\widetilde m_o(Z) + \widetilde m_i(X)- d_{\ol A_0}(Z,X)$.  Moreover,
since $D_0+D$ is $k$-connected, $D$ contains at least $k-w(B)$ arcs
covering any $D_0$-one-way bi-set.  Therefore $D$ contains at least
$\widetilde p_1({\cal F})$ arcs covering $\cal F$.  Since the members
of $\cal F$ are one-way with respect to $D\sp *[Z,X]$, the tail of
such an arc is not in $Z$ and its head is not in $X$.  Therefore the
total number $\gamma $ of arcs of $D$ is at least $\widetilde
p_1({\cal F}) + \widetilde m_o(Z) + \widetilde m_i(X)- d_{\ol
A_0}(Z,X)$, and \eqref{(felt.F0)} follows.

Sufficiency.  Assume that \eqref{(felt.F0)} holds.  Define

$$ \hbox{ $N_0\sp +(u):= \{u\} \cup \{v:uv\in A_0\}$, \ $N_0\sp -(u):=
\{u\} \cup \{v:vu\in A_0\}$.  }\ $$

\begin{claim} \label{momi} \begin{equation}k\leq m_o(v) + \delta _0(v) \leq n-1, \ \ \
k\leq m_i(v) + \varrho _0(v) \leq n-1 \label{(momi)} \end{equation} for every
node $v\in V$.  \end{claim}

\proof{Proof.}We prove only the first half of \eqref{(momi)} since the second
half follows analogously.  By choosing ${\cal F}=\emptyset $,
$Z=\{v\},$ and $X=V$, \eqref{(felt.F0)} gives rise to $m_o(v) +
\widetilde m_i(V) - d_{\ol A_0}(v,V) \leq \gamma $. Since $\widetilde
m_i(V)=\gamma $ and $d_{\ol A_0}(v,V) = n-1 -\delta _0(v)$, we obtain
that $m_o(v) +\delta _0(v)\leq n-1$.

Let $B=(V-v, V-N_0\sp +(v))$, ${\cal F}= \{B\}$, $Z=V-v$, and
$X=\emptyset $. Then $\widetilde p_1({\cal F})= k - \delta _0(v)$ and
$d_{\ol A_0}(Z,X)=0$ from which \eqref{(felt.F0)} implies that $k -
\delta _0(v) + \widetilde m_o(V-v) \leq \gamma $ from which $k\leq
m_o(v) + \delta _0(v)$, as required.  $\bullet$\endproof \medskip

Observe that ${\cal B}_0$ is a crossing family of bi-sets (for
details, see \cite{FrankJ31}).  Recall that the function $p_1$ on
${\cal B}_0$ was defined in the theorem by

\begin{equation}p_1(B):= k-w(B) \ (=k +\vert B_I\vert -\vert B_O\vert ). \end{equation}

\begin{figure}[t]
\centering
\begin{subfigure}[b]{0.33\textwidth}
  \centering
  \includegraphics[width=\linewidth]{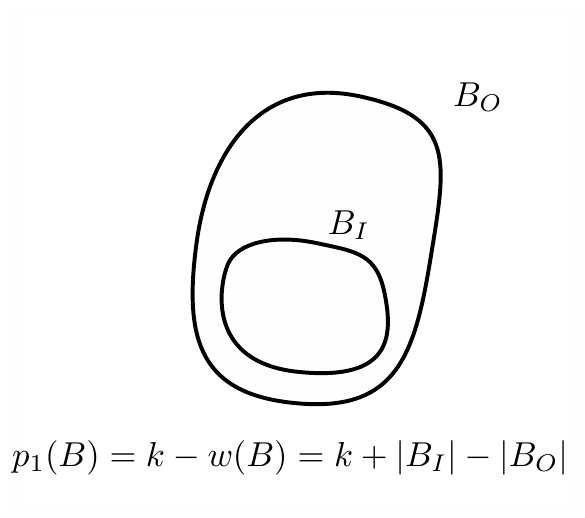}
  \caption{Definition of $p_1$}
  \label{fig:p1}
\end{subfigure}%
\begin{subfigure}[b]{.33\textwidth}
  \centering
  \includegraphics[width=\linewidth]{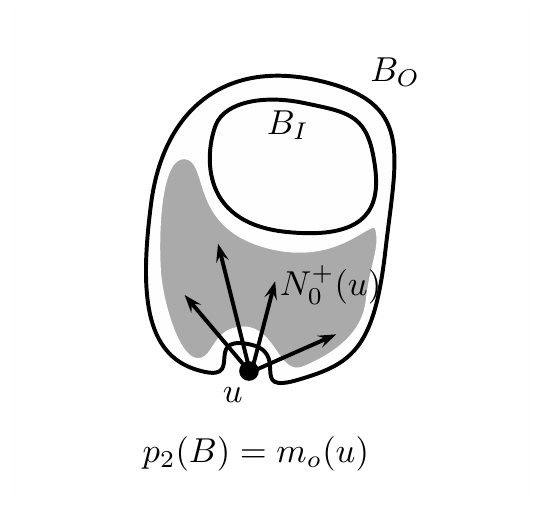}
  \caption{Definition of $p_2$}
  \label{fig:p2}
\end{subfigure}
\begin{subfigure}[b]{.33\textwidth}
  \centering
  \includegraphics[width=\linewidth]{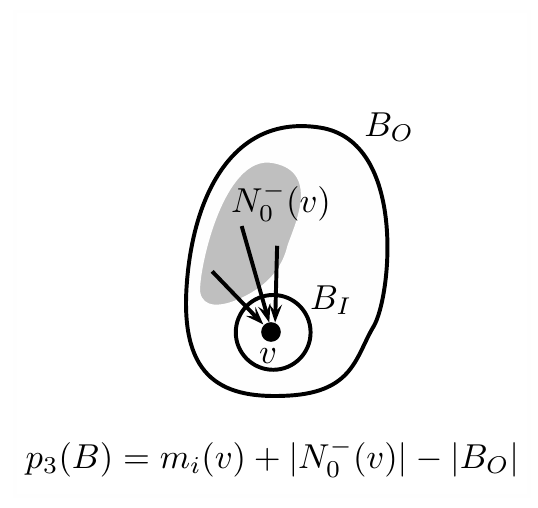}
  \caption{Definition of $p_3$}
  \label{fig:p3}
\end{subfigure}
\caption{Construction of the bi-set functions appearing in the definition of $p$}
\label{fig:p}
\end{figure}

In addition to $p_1$, we introduce three further functions defined on
${\cal B}_0$, as follows (see Figure~\ref{fig:p}).  For $B=(B_O,B_I)\in {\cal B}_0$, let

\begin{equation}p_2(B):= \begin{cases}
m_o(u) & \text{if $B_O=V-u, \ B_I=V-N_0\sp +(u)$ for
some $u\in V$}\\ 0 & \text{otherwise,}\end{cases} \end{equation}

\begin{equation} p_3(B):= \begin{cases}
 m_i(v)+\vert N_0\sp -(v)\vert - \vert B_O\vert &
\text{if $B_I=\{v\}$ for some $v\in V$}\\  0 & \text{otherwise,}\end{cases} \end{equation}

\begin{equation}p(B):= \max \{p_1(B),\ p_2(B),\ p_3(B), \ 0\}.  \end{equation}

Here $p_2$ and $p_3$ are to encode the out-degree and the in-degree
specifications, respectively.  Note, however, that the definitions of
$p_2$ and $p_3$ are not symmetric since a node $u$ determines a single
bi-set $B$ with outer set $V-u$ for which $p_2(B)$ is positive, while
a node $v$ may determine several bi-sets with inner set $\{v\}$ for
which $p_3(B)$ is positive.  The reason of this apparently undue
asymmetry is that both the supermodularity of $p$ and the simplicity
of the augmented digraph can be ensured only in this way.

\begin{claim} \label{p123} Let $B$ and $C$ be two crossing bi-sets.  Then
{\em (a)} \ $p_1(B)+p_1(C)= p_1(B\sqcap C)+ p_1(B\sqcup C)$ and {\em
(b)} \ $p_2(B)=0=p_2(C)$.  {\em (c)} \ If $p_3(B)>0$ and $p_3(C)>0$,
then $p_3(B)+p_3(C) = p_3(B\sqcap C)+ p_3(B\sqcup C)$.  \end{claim}

\proof{Proof.}Part (a) is immediate from the definition of $p_1$.  A bi-set
$B$ in ${\cal B}_0$ with $p_2(B)>0$ cannot cross any member of ${\cal
B}_0$ from which (b) follows.  In Part (c), $B_I=\{v\} =C_I$ for some
$v\in V$.  Therefore $B\sqcap C=(B_O\cap C_O,\{v\})$ and $B\sqcup
C=(B_O\cup C_O,\{v\})$ from which $p_3(B)+p_3(C) = p_3(B\sqcap C)+
p_3(B\sqcup C)$ follows.  $\bullet$\endproof

\begin{lemma} \label{supermodular} The bi-set function $p$ on ${\cal B}_0$ is
positively crossing supermodular.\end{lemma}

\proof{Proof.}Let $B$ and $C$ be two crossing members of ${\cal B}_0$ with
$p(B)>0, p(C)>0$.  By Part (b) of Claim \ref{p123}, $p(B)=\max
\{p_1(B),p_3(B)\}$ and $p(C)=\max \{p_1(C),p_3(C)\}$.  If
$p(B)=p_1(B)$ and $p(C)=p_1(C)$, then, by Part (a) of Claim
\ref{p123}, $p(B) + p(C)= p_1(B) + p_1(C)= p_1(B\sqcap C)+ p_1(B\sqcup
C) \leq p(B\sqcap C)+ p(B\sqcup C)$.  The supermodular inequality
follows analogously from Part (C) of Claim \ref{p123} when
$p(B)=p_3(B)$ and $p(C)=p_3(C)$.

Finally, suppose that $p(B) = p_3(B) =m_i(v)+\vert N_0\sp -(v)\vert -
\vert B_O\vert $ where $B_I=\{v\}$ for some $v\in V$, and $p(C)=p_1(C)
= k +\vert C_I\vert -\vert C_O\vert $. Now $C_I\cap B_I= \{v\}$ from
which $p_3(B\sqcap C)= m_i(v)+\vert N_0\sp -(v)\vert - \vert B_O\cap
C_O\vert $. Furthermore, $C_I\cup B_I= C_I$ from which $p_1(B\sqcup
C)= k +\vert C_I\vert - \vert B_O\cup C_O\vert $. Hence we have $$p(B)
+ p(C) = [m_i(v)+\vert N_0\sp -(v)\vert - \vert B_O\vert ] + [k +\vert
C_I\vert -\vert C_O\vert ] = $$ $$[m_i(v)+\vert N_0\sp -(v)\vert -
\vert B_O\cap C_O\vert ] + [k + \vert C_I\vert - \vert B_O\cup
C_O\vert ] \leq $$ $$p_3(B\sqcap C) + p_1(B\sqcup C)\leq p(B\sqcap C)
+ p(B\sqcup C). \qquad\bullet $$
\endproof

\begin{lemma} \label{pfed} A loopless digraph $D=(V,A)$ covering $p$ has at
least $\gamma $ arcs.  If $D$ has exactly $\gamma $ arcs, then {\em
(a)} \ \ $D$ fits $(m_o,m_i)$, {\em (b)} \ $D_0+D$ is $k$-connected,
and {\em (c)} \ $D_0+D$ is simple.\end{lemma}

\proof{Proof.}Since $D$ covers $p$ and $p\geq p_2$, we obtain that $D$ has at
least $m_o(u)$ arcs leaving $u$ for every node $u$, and therefore \begin{equation}
\vert A\vert = \sum [\delta _D(u):u\in V]= \sum [\varrho _D(V-u):u\in
V] \geq \widetilde m_o(V)=\gamma . \label{(atleast.gamma)} \end{equation}

Suppose now that $D$ has exactly $\gamma $ arcs.  Then
\eqref{(atleast.gamma)} implies $\delta _D(u)=m_o(u)$ for every node
$u\in V.$ Furthermore, $p_3(B)= m_i(v)$ for $B=(N_0\sp -(v), \{v\})$,
and thus $\varrho _D\geq p\geq p_3$ implies that $\varrho
_D(B)=\varrho _D(v)\geq m_i(v)$, from which $\varrho _D(v)=m_i(v)$ for
every $v\in V$.  That is, $D$ fits $(m_o,m_i)$.

By Lemma \ref{k-con}, in order to see that $H:=D_0+D$ is
$k$-connected, it suffices to show that $w(B)\geq k$ for every
non-trivial $H$-one-way bi-set $B$.  But this follows from $0=\varrho
_D(B)\geq p(B)\geq p_1(B) = k -w(B)$.

Finally, we prove that $D_0+D$ is simple.  Let $v$ be any node.  For
$B=(N_0\sp -(v),\{v\})$ we have $\varrho _D(v)\geq \varrho _D(B)\geq
p(B)\geq p_3(B)= m_i(v)=\varrho _D(v)$ from which one has equality
throughout.  But $\varrho _D(v)=\varrho _D(B)$ implies that every arc
$e=uv$ of $D$ must enter $N_0\sp -(v)$ as well, implying that $e$
cannot be parallel to any arc of $D_0$.  Suppose now that there is an
arc in $D$ from $u$ to $v$.  We have just seen that $e$ is not
parallel to any arc of $D_0$, that is, $u\not \in N_0\sp -(v)$.  Let
$B'=(N_0\sp -(v)+u, \{v\})$.  Then $p_3(B')=m_i(v)-1$ and hence
$$\varrho _D(v)\geq \varrho _D(B')+1 \geq p(B')+1 \geq p_3(B')+1 =
m_i(v) -1 +1 =\varrho _D(v)$$ from which one has equality throughout.
In particular, $\varrho _D(v)=\varrho _D(B')+1$, implying that there
is at most one arc in $D$ from $u$ to $v$.  $\bullet$\endproof

\begin{lemma} \label{novel3.mainlemma} There is a loopless digraph $D=(V,A)$
with $\gamma $ arcs covering $p$.  \end{lemma}

\proof{Proof.}Suppose indirectly that no such a digraph exists.  Theorem
\ref{arc-covering} implies that there is an independent family ${\cal
I}$ of non-trivial $D_0$-one-way bi-sets for which $\widetilde p({\cal
I})>\gamma $. We may assume that $\widetilde p({\cal I})$ is as large
as possible, modulo this, $\vert {\cal I}\vert $ is minimal, and
within this, \begin{equation}\hbox{ $ \sum [ \vert B_I\vert :  (B_O,B_I)\in {\cal
I}]$ is as small as possible.}\ \label{(kisbelso.3)} \end{equation} The
minimality of $\cal I$ implies $p(B)>0$ for every $B\in {\cal I}.$

\begin{claim} \label{kinai.egybeta.3} There are no two members
$B=(B_O,\{v\})$ and $C=(C_O,\{v\})$ of $\cal I$ (with the same inner
set $\{v\}$) for which $p(B)=p_3(B)$ and $p(C)=p_3(C)$.\end{claim}

\proof{Proof.}On the contrary, suppose that there are two such members.  Then
$B_O\cup C_O=V$ by the independence of $\cal I$.  If we replace the
two members $B$ and $C$ of $\cal I$ by the single bi-set $B\sqcap
C=(B_O\cap C_O,\{v\})$, then the resulting family ${\cal I}'$ is also
independent since any arc covering $B\sqcap C$ covers at least one of
$B$ and $C$.

Recall that $p_3(B\sqcap C)=m_i(v) + \vert N_0\sp -(v)\vert -\vert
B_O\cap C_O\vert $. By Claim \ref{momi}, $m_i(v)+ \varrho _0(v) \leq
n-1$ and hence $$p(B)+p(C)= p_3(B)+p_3(C) = m_i(v)+\vert N_0\sp
-(v)\vert -\vert B_O\vert + m_i(v)+\vert N_0\sp -(v)\vert -\vert
C_O\vert =$$ $$m_i(v)+\vert N_0\sp -(v)\vert -\vert B_O\cap C_O\vert
+m_i(v)+\vert N_0\sp -(v)\vert -\vert B_O\cup C_O\vert = $$ $$
p_3(B\sqcap C) + m_i(v) + \varrho _0(v)+1 - \vert V\vert \leq
p_3(B\sqcap C) +0 \leq p(B\sqcap C).$$

By the maximality of $\widetilde p({\cal I})$ we must have
$p(B)+p(C)=p(B\sqcap C)$ and hence $\widetilde p({\cal I}) =\widetilde
p({\cal I}')$ but this contradicts the minimality of $\vert {\cal
I}\vert $. $\bullet$\endproof \medskip

Let
\begin{eqnarray*}
&{\cal I}_1:=\{B\in {\cal I}:  p(B)=p_1(B) > \max
\{p_2(B),p_3(B)\} \},&\\[8pt]
&{\cal I}_2:=\{B\in {\cal I}:  p(B)=p_2(B)\},&\\[8pt]
&{\cal I}_3:=\{B\in {\cal I}:  p(B)=p_3(B)>p_2(B)\}.&
\end{eqnarray*}

Note that if $p_2(B)=p_3(B)\geq p_1(B)$ for $B\in {\cal I}$, then
$B\in {\cal I}_2$.  It follows that ${\cal I}_1, {\cal I}_2$, and
${\cal I}_3$ form a partition of $\cal I$.

Let $Z$ consist of those nodes $u$ for which the bi-set $(V-u,V-N_0\sp
+(u))$ is in ${\cal I}_2$.  Let $X$ consist of those nodes $v$ for
which there is a set $K_v$ such that $(K_v,\{v\})\in {\cal I}_3$.  By
Claim \ref{kinai.egybeta.3}, there is at most one such $K_v$ for each
$v\in X$.

\begin{claim} \label{kinai.Kv.3} For $v\in X$, one has $K_v=Z\cup N_0\sp
-(v)$.  \end{claim}

\proof{Proof.}$N_0\sp -(v) \subseteq K_v$ holds since $(K_v,\{v\})\in {\cal
I}_3\subseteq {\cal B}_0$.  If an element $u\in Z-N_0\sp -(v)$ would
not be in $K_v$, then $uv$ would cover both $(V-u,V-N_0\sp +(u))\in
{\cal I}_2$ and $(K_v,\{v\})\in {\cal I}_3$, contradicting the
independence of $\cal I$.  Therefore $Z\cup N_0\sp -(v)\subseteq K_v$.

To see the reverse inclusion suppose indirectly that there is an
element $u\in K_v-(Z\cup N_0\sp -(v))$.  Replace the member
$B=(K_v,v)$ of $\cal I$ by $B'=(K_v-u,v)$.  Note that $B'\in {\cal
B}_0$.  Since $p(B')\geq p_3(B') = p_3(B)+1 = p(B)+1$, the maximality
of $\widetilde p({\cal I})$ implies that the resulting system ${\cal
I}'$ is not independent.  Therefore there is a member $C$ of $\cal I$
that is covered by arc $uv$.

Since $u\not \in Z$, the bi-set $C$ cannot be in ${\cal I}_2$.  Since
$C$ and $B$ are distinct, Claim \ref{kinai.egybeta.3} implies that $C$
cannot be in ${\cal I}_3$ either.  Therefore $C\in {\cal I}_1$.  We
claim that $\vert C_I\vert \geq 2$.  Indeed, if $\vert C_I\vert =1$,
then $C_I=\{v\}$.  Since $p_1(C)>p_3(C)$, $p_1(C) = k-\vert C_O\vert +
\vert C_I\vert = k-\vert C_O\vert + 1$, and $p_3(C)=m_i(v) -\vert
C_O\vert + \vert N_0\sp +(v)\vert $, we obtain that $k + 1 > m_i(v) +
\vert N_0\sp +(v)\vert = m_i(v) + \varrho _0(v) +1$, contradicting the
second half of \eqref{(momi)}.

By replacing the member $C$ of ${\cal I}'$ with $C':=(C_O,C_I-v)$, we
obtain an independent family ${\cal I}''$.  Note that $C'\in {\cal
B}_0$.  Since $p(C')\geq p_1(C')=p_1(C)-1=p(C)-1$, we must have
$\widetilde p({\cal I}'')= \widetilde p({\cal I})$, but this
contradicts the minimality property given in \eqref{(kisbelso.3)}.
$\bullet$\endproof

\begin{claim} \label{p3.I3} $\widetilde p_3({\cal I}_3) = \widetilde m_i(X) -
d_{\ol A_0}(Z,X)$.  \end{claim}

\proof{Proof.}Let $B=(K_v,\{v\})\in {\cal I}_3$.  By Claim \ref{kinai.Kv.3},
$K_v=Z\cup N_0\sp -(v)$ and hence $p_3(B)=m_i(v)+\vert N_0\sp
-(v)\vert - \vert B_O\vert = m_i(v) - \vert Z-N_0\sp -(v)\vert =m_i(v)
-d_{\ol A_0}(Z,\{v\})$, and therefore $\widetilde p_3({\cal I}_3) =
\widetilde m_i(X) - d_{\ol A_0}(Z,X)$.  $\bullet$\endproof

\medskip By the definition of sets $Z$ and $X$ and by Claim
\ref{kinai.Kv.3}, we have $${\cal I}_2 = \{(V-u,V-N_0\sp +(u)):  \
u\in Z\} \ \hbox{ and }\ \ {\cal I}_3 = \{(Z\cup N_0\sp -(v), \{v\}):
\ v\in X\}.$$

\begin{claim} \label{one-way} Each member $B_1=(B_O,B_I)$ of~ ${\cal I}_1$ is
a one-way bi-set with respect to $D_0 + D\sp *[Z,X]$.\end{claim}

\proof{Proof.}Clearly, $B_1$ is $D_0$-one-way since ${\cal I}_1\subseteq
{\cal B}_0$.  Suppose now, indirectly, that an arc $e=uv$ of $D\sp
*[Z,X]$, which is not in $D_0$, covers $B_1$.  If $u\in Z$, then the
bi-set $B_2:=(V-u,V-N_0\sp +(u))$ belongs to ${\cal I}_2$ and $e$
covers $B_2$, contradicting the independence of $\cal I$ (namely, the
independence of $B_1$ and $B_2$).  If $v\in X$, then the bi-set
$B_3:=(Z\cup N_0\sp -(v), \{v\})$ belongs to ${\cal I}_3$ and $e$
covers $B_3$, contradicting the independence of $\cal I$ (namely, the
independence of $B_1$ and $B_3$).  $\bullet$\endproof \medskip

Since $\widetilde p_2({\cal I}_2)= \widetilde m_o(Z)$, we can conclude
that $$\gamma <\widetilde p({\cal I})= \widetilde p({\cal I}_1) +
\widetilde p({\cal I}_2) + \widetilde p({\cal I}_3) =$$ $$\widetilde
p_1({\cal I}_1) + \widetilde p_2({\cal I}_2) + \widetilde p_3({\cal
I}_3) = \widetilde p_1({\cal I}_1) + \widetilde m_o(Z) + \widetilde
m_i(X) - d_{\ol A_0}(Z,X),$$ and this inequality contradicts the
hypothesis \eqref{(felt.F0)} of the theorem since Claim \ref{one-way}
implies that ${\cal F}:={\cal I}_1$ consists of bi-sets which are
one-way with respect to $D_0+D\sp *[Z,X]$.  $\bullet$ $\bullet$\endproof

\medskip

\begin{remark} With a similar technique, it is possible to solve the
degree-specified node-connectivity augmentation problem when only the
augmenting digraph is required to be simple (and not the whole
augmented digraph).  Even more, one may prescribe a subset
$F_0\subseteq A_0$ of arcs of the starting digraph $D_0$ and request
for the degree-specified augmenting digraph $D$ we are looking for to
be simple and have no parallel arcs with the elements of $F_0$.  If
$F_0$ is empty, then this requires that the augmenting digraph be
simple, while if $F_0$ is the whole $A_0$, then this requires that the
augmented digraph be simple.  This can be done by revising first the
original definition of $N_0\sp +(u)$ and $N_0\sp -(u)$ as follows:  $$
\hbox{ $N_0\sp +(u):= \{u\} \cup \{v:uv\in F_0\}$ \ and \ $N_0\sp
-(u):= \{u\} \cup \{v:vu\in F_0\}$.  }\ $$ and applying then these
revised functions in defining $p_2$ and $p_3$.  The details are left
to the reader.  \end{remark}

\medskip

\begin{remark} In the derivation of Theorem
\ref{novel3.main} above, we applied the supermodular
arc-covering theorem for the general case when the bi-set function in
question is positively crossing supermodular.  But it is possible
to modify the definition of the bi-set function to be crossing
supermodular.  This is important in order to construct an algorithm.
The details shall be worked out in \cite{Berczi-Frank17}. \end{remark}

\medskip

\subsection{Simplified characterization for \texorpdfstring{$k=1$}{k=1}}

In Theorem \ref{erosnov} we considered the augmentation problem in
which an initial digraph $D_0$ was to be made strongly connected by
adding a degree-specified simple digraph $D$.  In Theorem
\ref{novel3.main}, the general degree-specified node-connectivity
augmentation problem was solved when not only the augmenting but the
augment{\bf ed} digraph was required to be simple.  The
characterization in Theorem \ref{novel3.main} had, however, an
aesthetic drawback in the sense that it included independent families
of bi-sets.  The goal of the present section is to show that in the
special case of $k=1$ this drawback can be eliminated.

\begin{theorem} \label{strong.main} Let $(m_o,m_i)$ be a degree-specification
with $\widetilde m_o(V)=\widetilde m_i(V)=\gamma $ and let
$D_0=(V,A_0)$ be a simple digraph with in- and out-degree functions
$\varrho _0$ and $\delta _0$, respectively.  There is a digraph
$D=(V,A)$ fitting $(m_o,m_i)$ for which $D_0\sp +:=D_0+D$ is simple
and strongly connected if and only if \begin{equation}\widetilde m_o(Z) + \widetilde
m_i(X)- d_{\ol A_0}(Z,X) \leq \gamma \ \hbox{holds for every
$X,Z\subseteq V$,}\ \label{(strong.felt.1)} \end{equation} \begin{equation}\hbox{
$\widetilde m_i(K)\geq 1$ \ and \ $\widetilde m_o(V-K)\geq 1$ \
whenever $\varrho _0(K)=0, \ \emptyset \subset K\subset V$, }\
\label{(strong.felt.12)} \end{equation} and \begin{equation}\widetilde m_o(Z) + \widetilde
m_i(X)- d_{\ol A_0}(Z,X) +1 \leq \gamma \label{(strong.felt.2)} \end{equation}
holds for every pair of disjoint non-empty subsets $X,Z\subset V$ for
which there is no dipath in $D_0$ from $X$ to $Z$, where $d_{\ol
A_0}(Z,X)$ denotes the number of arcs $a=zx\in \ol A_0$ for which
$z\in Z$ and $x\in X$.  \end{theorem}

\proof{Proof.}Necessity.  Let $D$ be a requested digraph.  Since $D_0+D$ is
simple, $D$ must be a subgraph of $\ol D_0$, and Condition
\eqref{(strong.felt.1)} is a special case of \eqref{(Ore)} when $H=\ol
D_0$.

To see the necessity of \eqref{(strong.felt.12)}, observe that, as
$D+D_0$ is strong, $D_0$ must admit an arc entering $K$ from which
$\widetilde m_i(K)\geq 1$ and $\widetilde m_o(V-K)\geq 1$.

Finally, consider the necessity of \eqref{(strong.felt.2)}.  Since $Z$
is not reachable from $X$ in $D_0$, there is a subset $K$ for which
$Z\subseteq K\subseteq V-X$ and $\varrho _0(K)=0$.  The digraph $D$
must have at least $\widetilde m_o(Z) + \widetilde m_i(X)- d_{\ol
A_0}(Z,X)$ arcs having tail in $Z$ or having head in $X$, and $D$ has
at least one more arc entering $K$.

Sufficiency.  Let ${\cal K}=\{K:  \ \emptyset \subset K\subset V,
\varrho _0(K)=0\}$.  Theorem \ref{novel3.main}, when applied in the
special case $k=1$, states that the requested digraph $D$ exists if
and only if \begin{equation}\widetilde p_1({\cal F}) + \widetilde m_o(Z) +
\widetilde m_i(X)- d_{\ol A_0}(Z,X) \leq \gamma
\label{(strong.felt.F0)} \end{equation} for subsets $X,Z\subseteq V$ and for
independent families $\cal F$ of non-trivial bi-sets which are one-way
with respect to $D_0+D\sp *[Z,X]$, where $p_1(B_O,B_I)=1-w(B)$ for
$B\in {\cal F}$.

Now $p_1(B_O,B_I)$ can be positive only if $B_O=B_I$.  The requirement
that a $(B_O,B_I)$ is $D_0$-one-way is equivalent to require that no
arc of $D_0$ enters $B_I$.  Therefore the condition can be restated as
follows.

\begin{equation}\vert {\cal F}\vert + \widetilde m_o(Z) + \widetilde m_i(X)-
d_{\ol A_0}(Z,X) \leq \gamma \label{(strong.felt.F01)} \end{equation} where
$\cal F$ is an independent family of sets $K\in {\cal K}$ such that no
arc of $D\sp *[Z,X]$ enters $K$.  This last property requires that
$Z\subseteq K$ and $X\subseteq V-K$.  The independence of $\cal F$
means that $\cal F$ consists of pairwise disjoint or pairwise
co-disjoint sets.

Our goal is to prove that \eqref{(strong.felt.F01)} follows from
\eqref{(strong.felt.1)}, \eqref{(strong.felt.12)}, and
\eqref{(strong.felt.2)}.  When $\cal F$ is empty,
\eqref{(strong.felt.F01)} is just \eqref{(strong.felt.1)}, and hence we
can assume that $\cal F$ is non-empty.  If neither $Z$ nor $X$ is
empty, then $\cal F$ has exactly one member, and in this case
\eqref{(strong.felt.F01)} and \eqref{(strong.felt.2)} coincide.  Hence
we can assume that at least one of $X$ and $Z$ is empty.  In analyzing
these cases, we rely on Condition \eqref{(strong.felt.12)} requiring
that $\widetilde m_i(K)\geq 1$ and $\widetilde m_o(V-K)\geq 1$ for
each $K\in {\cal K}$.

Suppose first that $Z=\emptyset $ and $X\not =\emptyset $. Then $\cal
F$ is a subpartition and $$\vert {\cal F}\vert + \widetilde m_o(Z) +
\widetilde m_i(X)- d_{\ol A_0}(Z,X) = \vert {\cal F}\vert + \widetilde
m_i(X) \leq \sum [\widetilde m_i(K):  K\in {\cal F}] + \widetilde
m_i(X) \leq \widetilde m_i(V)=\gamma ,$$ that is,
\eqref{(strong.felt.F01)} holds.

Suppose now that $Z\not =\emptyset $ and $X=\emptyset $. Then the
members of $\cal F$ include $Z$ and are pairwise co-disjoint.  Hence
$$\vert {\cal F}\vert + \widetilde m_o(Z) + \widetilde m_i(X)- d_{\ol
A_0}(Z,X) = \vert {\cal F}\vert + \widetilde m_o(Z) \leq \sum
[\widetilde m_o(V-K):  K\in {\cal F}] + \widetilde m_o(Z) \leq
\widetilde m_o(V)=\gamma ,$$ that is, \eqref{(strong.felt.F01)} holds.

Finally, suppose that $X=\emptyset =Z$.  When $\cal F$ is a
subpartition, we have $\vert {\cal F}\vert \leq \sum [\widetilde
m_i(K):  K\in {\cal F}] \leq \widetilde m_i(V)=\gamma $. When the
members of $\cal F$ are pairwise co-disjoint, we have $\vert {\cal
F}\vert \leq \sum [\widetilde m_i(V-K):  K\in {\cal F}] \leq
\widetilde m_i(V)=\gamma .$ Therefore in this case
\eqref{(strong.felt.F01)} holds.

Theorem \ref{novel3.main} ensures that the requested digraph $D$ does
exist.  $\bullet$ $\bullet$\endproof

\medskip

\begin{remark} Theorem \ref{strong.main} actually implies Theorem
\ref{erosnov} as follows.  Subdivide each arc of $D_0$ by a node and
define $m_o(z)=0$ and $m_i(z)=0$ for each subdividing node, and apply
Theorem \ref{strong.main} to the resulting digraph $D_0'$. \end{remark}

\medskip

\section{Degree-sequences of simple \texorpdfstring{$k$}{k}-connected digraphs}

In 1972, Wang and Kleitman \cite{Wang-Kleitman} characterized the
degree-sequences of simple $k$-connected undirected graphs.  The goal
of this section is to solve the analogous problem for directed graphs.
Before formulating the main result, we start with some preparatory
work.  Assume throughout that $1\leq k<n=\vert V\vert $. A node $v_f$
of a digraph $H$ is said to be {\bf full} if both $v_fu$ and $uv_f$
are arcs of $H$ for every node $u\not =v_f$.

\subsection{Preparations} \label{ulemma}

Let $Z$ and $X$ be two proper (but possibly empty) subsets of $V.$ Let
$D_{1}=(V, A_{1})$ denote the simple digraph in which $uv$ \ ($u\not
=v$) is an arc if $u\in Z$ or $v\in X$, that is, $D_1$ is the same as
$D\sp *[Z,X]$.  Note that each node in $X\cap Z$ is full.  Let ${\cal
B}_1$ denote the set of non-trivial $D_1$-one-way bi-sets.  Clearly,
$Z\subseteq B_O$ and $B_I\cap X=\emptyset $ hold for each $B\in {\cal
B}_1$.

Let ${\cal F}\subseteq {\cal B}_1$ be an independent family meaning
that each arc of the complete digraph $D\sp *=(V,A\sp *)$ covers at
most one member of $\cal F$.  Let $p_1(B)=k-w(B)$ where $w(B)=\vert
B_O-B_I\vert $. In proving our characterization of degree-sequences
realizable by simple $k$-connected digraphs, we shall need a simple
upper bound for $\widetilde p_1({\cal F})$.  Note that if $A$ is a set
of arcs for which $D_1\sp +=(V,A_1+A)$ is $k$-connected, then
$\widetilde p_1({\cal F})\leq \vert A\vert $.

\begin{lemma} \label{XZnagy} For an independent family ${\cal F}\subseteq
{\cal B}_1$, \begin{equation}\hbox{$\widetilde p_1({\cal F})\leq 0$ when $\vert
X\cap Z\vert \geq k,$}\ \end{equation} \begin{equation}\hbox{$\widetilde p_1({\cal F})\leq
k-\vert X\cap Z\vert $ when $\vert X\cap Z\vert <k$ and $\vert X\vert
,\vert Z\vert \geq k$.}\ \end{equation} \end{lemma}

\proof{Proof.}If $\vert X\cap Z\vert \geq k$, then $D_1$ has $k$ full nodes
and hence $D_1$ is $k$-connected, implying that $\widetilde p_1({\cal
F})\leq 0$.

The second part follows once we show that $D_1$ can be made
$k$-connected by adding $k'=k-\vert X\cap Z\vert $ new arcs.  By
symmetry, we may assume that $\vert Z\vert \geq \vert X\vert .$ Let
$x_1,\dots ,x_{k'}$ be distinct nodes in $X-Z$ and let $z_1,\dots
,z_{k'}$ be distinct nodes in $Z-X$.  Let $A=\{x_1z_1,\dots
,x_{k'}z_{k'}\}$ be a set of $k'$ disjoint arcs.

We claim that $D\sp +=(V,A_1+A)$ is $k$-connected.  Indeed, if $D\sp
+-K$ is not strong for some $K\subset V$, then $K$ contains every full
node and hence $X\cap Z\subseteq K$.  Moreover, $K$ must hit every arc
in $A$ since if $\{x_i,z_i\}\cap K=\emptyset $ for some $i$, then
$D\sp +-K$ would be strong as there is an arc $z_iu$ and an arc $ux_i$
for each node $u\in V-K-\{x_i,z_i\}$.  Therefore $\vert K\vert \geq k$
and hence $D\sp +$ is $k$-connected.  This means that $D_1$ has been
made $k$-connected by adding $k'$ new arcs, from which $\widetilde
p_1({\cal F})\leq k'= k-\vert X\cap Z\vert $. $\bullet$\endproof \medskip

The total out-deficiency of the nodes is defined by $\sigma _o= \sum
[(k-\delta _{D_{1}}(v))\sp +:v\in V]$.  Clearly, if $\vert X\vert <k$,
then $$ \sigma _o= (n-\vert Z\vert )(k-\vert X\vert ) +\vert X-Z\vert
.$$

\begin{lemma} \label{XZkicsi} Let ${\cal F}\subseteq {\cal B}_1$ be an
independent family.  If $\vert X\vert <k$ and $\vert Z\vert \geq \vert
X\vert $, then \begin{equation}\widetilde p_1({\cal F})\leq \sigma _o.
\label{(Xkicsi)} \end{equation} \end{lemma}

\proof{Proof.}The number of arcs of $D_1$ can be expressed as follows.  \begin{equation}
\vert A_1\vert = \vert Z\vert (n-1)+(n-\vert Z\vert )\vert X\vert
-\vert X-Z\vert . \label{(A1.elszam)} \end{equation}

The inequality \eqref{(Xkicsi)} is equivalent to \begin{equation}\sum [k - w(B):
B\in {\cal F}] \leq (n-\vert Z\vert )(k-\vert X\vert ) +\vert X-Z\vert
. \label{(Xkicsi.ekv)} \end{equation}

Let $q:=\vert \cal F\vert $. We distinguish two cases.  \medskip

{\bf \noindent Case 1} \ $q\geq n-\vert Z\vert $.

\begin{claim} \label{1.eset} \begin{equation}\sum [(n-1) - w(B):  B\in {\cal F}] \leq
(n-\vert Z\vert )(n-1-\vert X\vert ) +\vert X-Z\vert .
\label{(1.eset)} \end{equation} \end{claim}

\proof{Proof.}Since $\cal F$ is independent, the total number of those arcs
of the complete digraph $D\sp *=(V,A\sp *)$ which cover a member $B$
of $\cal F$ is $\sum [\vert B_I\vert (n-\vert B_O\vert ):  B\in {\cal
F}]$.  Since $A_1$ covers no member of $\cal F$, we conclude that
$$\sum [\vert B_I\vert (n-\vert B_O\vert ) :  B\in {\cal F}] + \vert
Z\vert (n-1)+(n-\vert Z\vert )\vert X\vert -\vert X-Z\vert =$$ $$ \sum
[\vert B_I\vert (n-\vert B_O\vert ) :  B\in {\cal F}] + \vert A_1\vert
\leq \vert A\sp *\vert = n(n-1).$$

By observing that $$\vert B_I\vert (n-\vert B_O\vert )\geq 1 \cdot
[n-\vert B_O\vert +(\vert B_I\vert -1)] = n-1 - w(B),$$ we obtain $$
\sum [(n-1) - w(B):  B\in {\cal F}] \leq \sum [\vert B_I\vert (n-\vert
B_O\vert ):  B\in {\cal F}] \leq $$ $$ n(n-1) - [ \vert Z\vert
(n-1)+(n-\vert Z\vert )\vert X\vert -\vert X-Z\vert ] = (n-\vert
Z\vert )(n-1-\vert X\vert ) +\vert X-Z\vert , $$ as required for
\eqref{(1.eset)}.  $\bullet$\endproof \medskip

As we are in Case 1, $q(n-1-k)\geq (n-\vert Z\vert )(n-1-k)$.  By
subtracting this inequality from \eqref{(1.eset)}, we obtain
\eqref{(Xkicsi.ekv)}, proving the lemma in Case 1. \medskip

{\bf \noindent Case 2} \ $q<n-\vert Z\vert $. Recall that we have
assumed $\vert Z\vert \geq \vert X\vert $ and $\vert X\vert <k$.

\begin{claim} \label{case2} \begin{equation}\sum [n-w(B) :  B\in {\cal F}]\leq q(n-\vert
X\vert ) +\vert X-Z\vert . \label{(case2)} \end{equation} \end{claim}

\proof{Proof.}Let $h:=\vert X-Z\vert $. Let $x_1,\dots ,x_{h}$ be the
elements of $X-Z$ and let $\{z_1,\dots ,z_{h}\}$ be a subset of $Z-X$.
Consider the set $A=\{x_1z_1,\dots ,x_{h}z_{h}\}$ of disjoint arcs.
For $B\in {\cal F}$, let $\alpha (B)$ denote the number of arcs in $A$
entering $B_O$ but not $B_I$, and $\beta (B)$ the number of arcs in
$A$ entering both $B_O$ and $B_I$ (that is, covering $B$).  Note that
$X\cap Z\subseteq B_O-B_I$.

Clearly, $\alpha (B)+\beta (B) = \vert X-B_O\vert $. The definition of
$\alpha (B)$ immediately shows that $\alpha (B)\leq w(B)-\vert X\cap
B_O\vert $. It follows that $$\beta (B) = \vert X-B_O\vert -\alpha (B)
\geq \vert X-B_O\vert - [w(B) -\vert X\cap B_O\vert ]=\vert X\vert -
w(B).$$

Since $\cal F$ is independent, every arc in $A$ covers at most one
member of $\cal F$ and hence $$\sum [\beta (B):  B\in {\cal F}] \leq
\vert A\vert =\vert X-Z\vert .$$ By combining these observations, we
obtain $$\vert X-Z\vert \geq \sum [\beta (B):  B\in {\cal F}] \geq
\sum [\vert X\vert -w(B):  B\in {\cal F}] = \sum [ n-w(B):  B\in {\cal
F}] - q(n-\vert X\vert ),$$ that is, $\sum [ n-w(B):  B\in {\cal
F}]\leq q(n-\vert X\vert ) +\vert X-Z\vert $, and the claim follows.
$\bullet$\endproof \medskip

By subtracting first the equality $\sum [n-k:  B\in {\cal F}] =
q(n-k)$ from \eqref{(case2)} and applying then the assumption
$q<n-\vert Z\vert $, one obtains $$\sum [k-w(B):  B\in {\cal F}]\leq
q(k-\vert X\vert ) +\vert X-Z\vert < (n-\vert Z\vert )(k-\vert X\vert
) +\vert X-Z\vert ,$$ that is, \eqref{(Xkicsi.ekv)} holds, completing
the proof of the lemma.  $\bullet$ $\bullet$\endproof

\subsection{The characterization}

\begin{theorem} \label{novel4.main} Let $V$ be a set of $n$ nodes and let
$1\leq k\leq n-1$.  Suppose that $(m_o,m_i)$ is a degree-specification
meeting \eqref{(kinai.alap)}.  There exists a simple $k$-connected
digraph fitting $(m_o,m_i)$ if and only if \begin{equation}\widetilde m_i(X) +
\widetilde m_o(Z) - \vert X\vert \vert Z\vert + k \leq \gamma \hbox{
whenever $X,Z\subset V$, $X\neq Z$. }\ \label{(k-con.m.1)} \end{equation} Moreover, it
suffices to require the inequality in \eqref{(k-con.m.1)} only for
its special case when $\vert X\cap Z\vert <k$, $X$ consists of the $h$
largest values of $m_i$ and $Z$ consists of the $j$ largest values of
$m_o$.  \end{theorem}

\proof{Proof.}Necessity.  Suppose that there is a digraph $D$ with the
requested properties, and let $X$ and $Z$ be two distinct, proper
subsets of $V$.  If $\vert X\cap Z\vert \geq k$, then
\eqref{(kinai.alap)} implies $\widetilde m_i(X) + \widetilde m_o(Z) -
\vert X\vert \vert Z\vert + k \leq \widetilde m_i(X) + \widetilde
m_o(Z) - \vert X\vert \vert Z\vert +\vert X\cap Z\vert \leq \gamma $
and \eqref{(k-con.m.1)} holds.  This argument also implies that, given
\eqref{(kinai.alap)}, it suffices to require the inequality in
\eqref{(k-con.m.1)} only for $X$ and $Z$ for which $\vert X\cap Z\vert
<k$.

Assume now that $\vert X\cap Z\vert \leq k-1$.  By the simplicity of
$D$, there are at most $\vert X\vert \vert Z\vert -\vert X\cap Z\vert
$ arcs with tail in $Z$ and head in $X$.  Therefore the total number
of arcs with tail in $Z$ or head in $X$ is at least $\widetilde m_i(X)
+\widetilde m_o(Z) - \vert X\vert \vert Z\vert +\vert X\cap Z\vert $.
Moreover, the $k$-connectivity of $D$ implies for distinct subsets
$X,Z\subset V$ that there are at least $k-\vert X\cap Z\vert $ arcs
from $V-Z$ to $Z-X$ or there are at least $k-\vert X\cap Z\vert $ arcs
from $X-Z$ to $V-X$.  Since the tails of these arcs are not in $Z$ and
the heads of these arcs are not in $X$, we can conclude that
$\widetilde m_i(X) + \widetilde m_o(Z) - \vert X\vert \vert Z\vert + k
= \widetilde m_i(X) + \widetilde m_o(Z) - \vert X\vert \vert Z\vert
+\vert X\cap Z\vert + (k - \vert X\cap Z\vert ) \leq \gamma $, that
is, the inequality in \eqref{(k-con.m.1)} holds in this case, too.

Sufficiency.  Suppose indirectly that the requested digraph does not
exist.  By applying Theorem \ref{novel3.main} to the empty digraph
$D_0=(V,\emptyset )$ and observing that in this case $d_{\ol
D_0}(Z,X)= \vert Z\vert \vert X\vert -\vert Z\cap X\vert $, we obtain
that there exists an independent family $\cal F$ of bi-sets and
subsets $Z$ and $X$ for which $Z\subseteq B_O, X\cap B_I=\emptyset ,$
and \begin{equation}\widetilde p_1({\cal F}) + \widetilde m_o(Z) + \widetilde
m_i(X) - \vert X\vert \vert Z\vert + \vert X\cap Z\vert >\gamma .
\label{(kinai.nagyobb.4b)} \end{equation}

If we reorient a digraph by reversing each of its arcs, then the
$k$-connectivity is preserved and the in-degrees and out-degrees
transform into each other.  Therefore the roles of $m_o$ and $m_i$ in
Theorem \ref{novel3.main} are symmetric and thus we may and shall
assume that $\vert Z\vert \geq \vert X\vert $ in
\eqref{(kinai.nagyobb.4b)}.

We can also assume that $\cal F$ is minimal and hence $p_1(B)>0$ for
each $B\in {\cal F}$.  Since $(m_o,m_i)$ is required to meet
\eqref{(kinai.alap)}, $\cal F$ is non-empty.  Lemma \ref{XZnagy}
implies $\vert X\cap Z\vert <k$.

Suppose first that $\vert X\vert ,\vert Z\vert \geq k$.  Since $\vert
X\cap Z\vert <k$ we have $X\not =Z$, $X\not =V$, and $Z\not =V$.  By
Lemma \ref{XZnagy}, $\widetilde p_1({\cal F})\leq k-\vert X\cap Z\vert
.$ From \eqref{(kinai.nagyobb.4b)}, we have

$$\gamma < \widetilde p_1({\cal F}) + \widetilde m_o(Z) + \widetilde
m_i(X) - \vert X\vert \vert Z\vert + \vert X\cap Z\vert \leq $$
$$k-\vert X\cap Z\vert + \widetilde m_o(Z) + \widetilde m_i(X) - \vert
X\vert \vert Z\vert + \vert X\cap Z\vert = k + \widetilde m_o(Z) +
\widetilde m_i(X) - \vert X\vert \vert Z\vert ,$$ contradicting
\eqref{(k-con.m.1)}.

In the remaining case $\vert Z\vert \geq \vert X\vert $ and $\vert
X\vert <k$.  Lemma \ref{XZkicsi} implies that $\widetilde p_1({\cal
F})\leq \sigma _o = (n-\vert Z\vert )(k-\vert X\vert ) +\vert X-Z\vert
.$ By applying the inequality in \eqref{(k-con.m.1)} to $X=\emptyset $
and $Z=V-v$, one gets $\widetilde m_o(v)\geq k$ from which $\widetilde
m_o(V-Z)\geq k\vert V-Z\vert $. By using again
\eqref{(kinai.nagyobb.4b)}, we obtain 
\begin{eqnarray*}
&\gamma < \widetilde p_1({\cal
F}) + \widetilde m_o(Z) + \widetilde m_i(X) - \vert X\vert \vert
Z\vert + \vert X\cap Z\vert \leq&\\[8pt]
&(k-\vert X\vert )(n-\vert
Z\vert ) +\vert X-Z\vert + \widetilde m_o(Z) + \widetilde m_i(X) -
\vert X\vert \vert Z\vert + \vert X\cap Z\vert =&\\[8pt]
&k(\vert V\vert
-\vert Z\vert ) -\vert X\vert \vert V\vert + \vert X\vert \vert Z\vert
+ \widetilde m_o(Z) + \widetilde m_i(X) - \vert X\vert \vert Z\vert +
\vert X\vert \leq&\\[8pt]
&\widetilde m_o(V-Z) -\vert X\vert \vert V\vert
+ \widetilde m_o(Z) + \widetilde m_i(X) + \vert X\vert =&\\[8pt]
&\widetilde m_o(V) +\widetilde m_i(X) -\vert X\vert \vert V\vert +
\vert X\cap V\vert ,&
\end{eqnarray*}
and this contradicts the inequality
\eqref{(kinai.alap)} with $V$ in place of $Z$.  $\bullet$\endproof

\medskip %{\tiny {\bf directory:  fenyo, file:  novel4, \today}}

\section{Degree-sequences of \texorpdfstring{$k$}{k}-elementary bipartite graphs}

As an application of Theorem \ref{novel4.main}, we extend a
characterization of degree-sequences of elementary bipartite graphs
(bigraphs, for short), due to R. Brualdi \cite{Brualdi80}, to
$k$-elementary bigraphs.  A simple bipartite graph $G=(S,T;E)$ is
called {\bf elementary} if $G$ is perfectly matchable and the union of
its perfect matchings is a connected subgraph.  It is known (see,
\cite{Lovasz-Plummer}, p. 122) that $G$ is elementary if and only if
it is either just one edge or $\vert S\vert =\vert T\vert \geq 2$ and
the Hall-condition holds with strict inequality for every non-empty
proper subset of $T$.  Another equivalent formulation requires that
either $G$ is just one edge or $G-s-t$ has a perfect matching for each
$s\in S$ and $t\in T$.  Yet another characterization states that $G$
is elementary if and only if it is connected and every edge belongs to
a perfect matching.  It should be noted that there is a one-to-one
correspondence between elementary bipartite graphs and fully
indecomposable $(0,1)$-matrices \cite{Brualdi-Ryser}.

Let $k\leq n-1$ be a positive integer.  We call a simple bigraph
$G=(S,T;E)$ with $\vert S\vert =\vert T\vert =n$ {\bf $k$-elementary}
if the removal of any $j$-element subset of $S$ and any $j$-element
subset of $T$ leaves a perfectly matchable graph for each $0\leq j\leq
k$.  Equivalently, $\vert \Gamma (X)\vert \geq \vert X\vert +k$ for
every subset $X\subseteq S$ with $\vert X\vert \leq n-k$, where
$\Gamma (X)$ denotes the set of neighbours of $X$, that is, $\Gamma
(X):=\{t:$ \ there is an edge $st\in E$ with $s\in X\}$.  Obviously,
this last property is equivalent to requiring $\vert \Gamma (X)\vert
\geq \vert X\vert +k$ for every subset $X\subseteq T$ with at most
$n-k$ elements.  Note that for $k=n-1$ the complete bigraph $K_{n,n}$
is the only $k$-elementary bigraph.  For $n\geq 2$, a bigraph is
$1$-elementary if and only if it is elementary, while for $n=1$ the
bigraph $K_{1,1}$ consisting of a single edge is considered elementary
but not 1-elementary.  The notion of $k$-elementary bigraphs was
introduced by Frank and V\'egh in \cite{FrankJ56} where they developed
a purely combinatorial algorithm the make a $(k-1)$-elementary bigraph
$k$-elementary by adding a minimum number of new edges.

\begin{figure}
\centering
\includegraphics[width=0.8\textwidth]{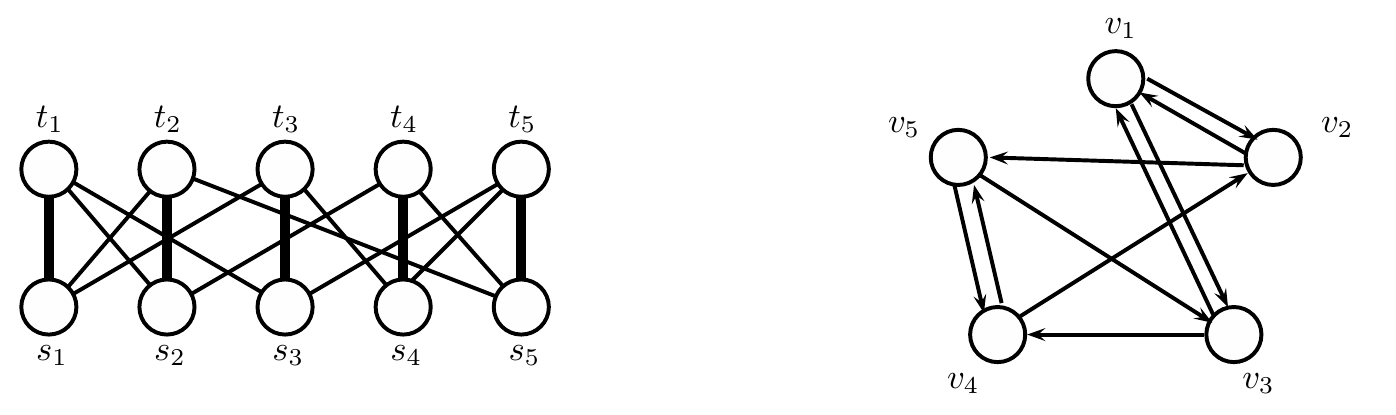}
\caption{A perfect matching of a 2-elementary bigraph and the associated digraph}
\label{fig:kelementary}
\end{figure}

There is a natural correspondence between $k$-connected digraphs and
$k$-elementary bigraphs with a specified perfect matching.  Let
$M=\{e_1,\dots ,e_n\}$ be a perfect matching of a $k$-elementary
bigraph $G=(S,T;E)$.  Suppose that the elements of $S$ and $T$ are
indexed in such a way that $e_i=s_it_i$ \ $(i=1,\dots ,n)$.  Let
$V=\{v_1,\dots ,v_n\}$ be a node-set where $v_i$ is associated with
$e_i$.  We associate a digraph $D=(V,A)$ with the pair $(G,M)$ as
follows.  For each edge $e=t_is_j$ of $G$ $(i\not =j)$ (that is, $e\in
E-M$) let $v_iv_j$ be an arc of $D$ (see Figure~\ref{fig:kelementary} for an example).  It can easily be observed that
the digraph $D$ obtained in this way is $k$-connected if and only if
$G$ is $k$-elementary.  By using network flow techniques, a digraph
can be checked in polynomial time whether it is $k$-connected or not,
and therefore a bipartite graph can also be checked for being
$k$-elementary.

Conversely, one can associate a perfectly matchable bigraph with a
digraph $D$ on node-set $\{v_1,\dots ,v_n\}$, as follows.  Define
$G=(S,T;E)$ so that $s_it_i$ belongs to $E$ for $i=1,\dots ,n$ and,
for every arc $v_jv_h$ of $D$, let $t_js_h$ be an edge of $G$.  This
$G$ is a simple $k$-elementary graph precisely if $D$ is
$k$-connected.

Let $m=(m_S,m_T)$ be a degree-specification for which $\widetilde
m_S(S)=\widetilde m_T(T)=\gamma $ and assume that the set ${\cal
G}(m_S,m_T)$ of simple bipartite graphs fitting $(m_S,m_T)$ is
non-empty, that is, by a theorem of Gale and Ryser \cite{Gale57},
\cite{Ryser57} \begin{equation}\hbox{ $\widetilde m_S(X) + \widetilde m_T(Z) -
\vert X\vert \vert Z\vert \leq \gamma $ whenever $X\subseteq S, \
Z\subseteq T$.}\ \label{(msmt.nemures)} \end{equation}

Brualdi \cite{Brualdi80} characterized the degree-sequences of
elementary bipartite graphs in terms of fully indecomposable
$(0,1)$-matrices.  Here we extend his results (apart from its trivial
special case when $n=1$) to $k$-elementary bipartite graphs.

\begin{theorem} \label{kelemi} Suppose that ${\cal G}(m_S,m_T)$ is non-empty,
that is, \eqref{(msmt.nemures)} holds.  There is a $k$-elementary
member $G$ of ${\cal G}(m_S,m_T)$ if and only if \begin{equation}\widetilde m_S(X)
+ \widetilde m_T(Z) - \vert X\vert \vert Z\vert + (n -\vert X\vert
-\vert Z\vert +k) \leq \gamma \ \hbox{whenever $X\subset S, \ Z\subset
T$.  }\ \label{(elemi.1)} \end{equation}

Furthermore, it suffices to require \eqref{(elemi.1)} only for
$X\subset S$ consisting of the $h$ largest $m_S$-valued elements and
for $Z\subset T$ consisting of the $j$ largest $m_T$-valued elements
$(h,j<n)$.  Moreover, if $m_S(s_1)\leq \cdots \leq m_S(s_n)$ and
$m_T(t_1)\geq \cdots \geq m_T(t_n)$, the graph $G$ can be chosen in
such a way that $s_1t_1,\dots ,s_nt_n$ is a perfect matching of $G$.
\end{theorem}

\proof{Proof.}Necessity.  Suppose that there is a requested bigraph $G$.
Since $G$ is $k$-elementary, the degree of each node is at least
$k+1$.  For $X\subset S$ and $Z\subset T$, let $\gamma _1$ denote the
number of edges incident to a node in $X\cup Z$ while $\gamma _2$ is
the number of the remaining edges.  Then $\gamma _1\geq \widetilde
m_S(X) + \widetilde m_T(Z) - \vert X\vert \vert Z\vert $.

If $\vert Z\vert \geq k$, then $T-Z$ has at least $\vert T-Z\vert +k$
neighbours and hence $T-Z$ has at least $\vert T-Z\vert +k-\vert
X\vert = n-\vert X\vert -\vert Z\vert +k$ neighbours in $S-X$, from
which $\gamma _2\geq \vert T-Z\vert +k-\vert X\vert $. If $\vert
Z\vert <k$, then each node in $S-X$ has at least $k+1-\vert Z\vert $
neighbours in $T-Z$ from which $\gamma _2\geq \vert S-X\vert
(k+1-\vert Z\vert )= (n-\vert X\vert )(k-\vert Z\vert ) + n-\vert
X\vert \geq k-\vert Z\vert +n-\vert X\vert $ in this case, too.
Therefore we have $$\gamma =\gamma _1+\gamma _2\geq \widetilde m_S(X)
+ \widetilde m_T(Z) - \vert X\vert \vert Z\vert + n-\vert X\vert
-\vert Z\vert +k,$$ and hence the inequality in \eqref{(elemi.1)}
holds.

Sufficiency.  We may suppose that the elements of $S$ and $T$ are
ordered in such a way that $m_S(s_1)\leq \cdots \leq m_S(s_n)$ and
$m_T(t_1)\geq \cdots \geq m_T(t_n)$.  The inequality in
\eqref{(elemi.1)}, when applied to $X:=S-s$ and to $Z=\emptyset $,
implies that $m_S(s)\geq k+1$ for $s\in S$ and $m_T(t)\geq k+1$
follows analogously for $t\in T$.  Let $V=\{v_1,\dots ,v_n\}$ be a
set.  Let $m_o(v_j):  = m_T(t_j)-1$ and $m_i(v_j):=m_S(s_j)-1$ for
$j=1,\dots ,n$.  Let $\gamma '=\gamma -n$.

Let $Z_j$ denote the first $j$ ($0\leq j\leq n$) elements of $V$ and
let $X_h$ denote the last $h$ ($0\leq h\leq n$) elements of $V$.  The
(possibly empty) subset of $T$ corresponding to $Z_j$ is denoted by
$Z$ while the subset of $S$ corresponding to $X_j$ is denoted by $X$.
Due to the assumption made on the ordering of the elements of $V$,
$Z_j$ consists of the $j$ elements of $V$ with largest $m_o$-values
while $X_h$ consists of the $h$ elements of $V$ with largest
$m_i$-values.

\begin{claim} \label{kelemi.1} Suppose that $j<n$ and $h<n$, that is,
$X\subset S$ and $Z\subset T$.  Then the inequality in
\eqref{(k-con.m.1)} with $\gamma '$ in place of $\gamma $ holds.
\end{claim}

\proof{Proof.}$\widetilde m_o(Z_j) + \widetilde m_i(X_h) -jh +k = \widetilde
m_T(Z) -\vert Z\vert + \widetilde m_S(X)-\vert X\vert -\vert X\vert
\vert Z\vert + k \leq \gamma - n = \gamma '.$ $\bullet$\endproof \medskip

\begin{claim} \label{kelemi.2} Suppose that $j\leq n$ and $h\leq n$, that is,
$X\subseteq S$ and $Z\subseteq T$.  Then the inequality in
\eqref{(kinai.alap)} with $\gamma '$ in place of $\gamma $ holds.
\end{claim}

\proof{Proof.}Suppose first that $h+j\geq n$, that is, $\vert X_h\cup
Z_j\vert =n$.  Then
\begin{eqnarray*}
&\widetilde m_i(X_h) + \widetilde m_o(Z_j) - \vert X_h\vert \vert
Z_j\vert + \vert X_h\cap Z_j\vert =&\\[8pt]
& \widetilde m_S(X) -\vert X\vert +
\widetilde m_T(Z)-\vert Z\vert -\vert X\vert \vert Z\vert + \vert
X_h\cap Z_j\vert =&\\[8pt]
& \widetilde m_S(X) + \widetilde m_T(Z) -\vert
X\vert \vert Z\vert -\vert X_h\vert -\vert Z_j\vert + \vert X_h\cap
Z_j\vert =&\\[8pt]
& \widetilde m_S(X) + \widetilde m_T(Z) -\vert X\vert
\vert Z\vert - \vert X_h\cup Z_j\vert\leq &\\[8pt]
& \gamma - \vert X_h\cup
Z_j\vert =\gamma -n=\gamma '.&
\end{eqnarray*}
Here the last inequality follows from
\eqref{(msmt.nemures)}.

Second, suppose that $h+j<n$, that is, $X_h\cup Z_j\subset V$ implying
that $X_h\cap Z_j=\emptyset $. By Claim \ref{kelemi.1}, $\widetilde
m_o(Z_j) + \widetilde m_i(X_h) -\vert X_h\vert \vert Z_j\vert +k \leq
\gamma ' $ and hence $\widetilde m_i(X_h) + \widetilde m_o(Z_j) -
\vert X_h\vert \vert Z_j\vert + \vert X_h\cap Z_j\vert = \widetilde
m_i(X_h) + \widetilde m_o(Z_j) - \vert X_h\vert \vert Z_j\vert <\gamma
'$.  $\bullet$\endproof

\medskip

Since the conditions \eqref{(kinai.alap)} and \eqref{(k-con.m.1)} were
shown above to hold, Theorem \ref{novel4.main} implies the existence
of a $k$-connected simple digraph $D=(V,A)$ fitting $(m_o,m_i)$.
Consider the bigraph $G=(S,T;E)$ associated with $D$ (in which
$s_it_i$ belongs to $E$ for $i=1,\dots ,n$ and, for every arc $v_jv_h$
of $D$, let $t_js_h$ be an edge of $G$).  This $G$ is a simple
$k$-elementary bigraph fitting $(m_S,m_T)$.  $\bullet$ $\bullet$\endproof

\medskip \medskip {\bf Acknowledgements} \ \ We are grateful to N\'ora
B\"odei for the valuable discussions on the topic and for her
indispensable help in exploring and understanding the paper
\cite{Hong-Liu-Lai}.  Special thanks are due to Zolt\'an Szigeti for
carefully checking the details of a first draft.  Thanks are also due
to Richard Brualdi for an ongoing email correspondence in which he
always promptly provided us with extremely useful background
information on the history of the topic.  The two anonymous referees
provided a great number of particularly useful comments.  We
gratefully acknowledge their invaluable efforts.

The research was supported by the Hungarian Scientific Research Fund -
OTKA, No K109240.  The work of the first author was financed by a
postdoctoral fellowship provided by the Hungarian Academy of Sciences.

\medskip %{\tiny {\bf directory:  fenyo, file:  kelemi, \today}}

\end{document}